\documentclass[11pt,oneside]{article}
\usepackage{amsfonts,amssymb,amsmath,verbatim,authblk,hyperref,graphicx,
amsbsy,mathrsfs,caption,subcaption,graphics,fancyhdr}

\newcommand{\beq}{\begin{equation}\ }
\newcommand{\eeq}{\end{equation}\ }
\newtheorem{lemma}{Lemma}[section]
\newtheorem{definition}{Definition}[section]
\newtheorem{corollary}{Corollary}[section]

\newcommand{\toprule}{\hline\noalign{\smallskip}}
\newcommand{\colrule}{\noalign{\smallskip}\hline\noalign{\smallskip}}
\newcommand{\botrule}{\noalign{\smallskip}\hline}

\begin{document}
\pagestyle{fancy}
\renewcommand\headrulewidth{0pt}
\lhead{}\chead{}\rhead{}
\cfoot{\vspace*{1.5\baselineskip}\thepage}

\author[1,2]{Demetris T. Christopoulos}
\affil[1]{\small{National and Kapodistrian University of Athens, Department of Economics}}
\affil[2]{\tt{dchristop@econ.uoa.gr}, \tt{dem.christop@gmail.com} }

\title{{\itshape Developing methods for identifying the inflection point of a convex/concave curve}}

\maketitle

\begin{abstract}
We are introducing two methods for revealing the true inflection point of data that contains or not error. The starting point is a set of geometrical properties that follow the existence of an inflection point p for a smooth function. These properties connect the concept of convexity/concavity before and after p respectively with three chords defined properly. Finally a set of experiments is presented for the class of sigmoid curves and for the third order polynomials.
\end{abstract}


\section{Introducing geometrical methods}
\label{intro}
The finding of the inflection point is performed by adopting a suitable model and then with regression or maximum likelihood estimation techniques. Another approach is that of the Differential Geometry view, when we define a proper measure of discrete curvature and then we choose the points with such a measure as close to zero can be, see for example \cite{mok-08}, \cite{han-01} or Gaussian smoothing techniques, see \cite{mok-86}. A review of the relevant shape analysis methods can be found in \cite{lon-98}.\\  
\par We are going to present two new methods for identifying the inflection point of any given convex/concave curve based on the definition and its geometrical properties only and without any regression or splines representation. We will use a \emph{generalization of bisection method in root finding} by choosing two points where the inflection point is between and then by taking the middle point as an estimator. The methods that will be presented will be able to iteratively converge to the actual inflection point in a manner similar to that of bisection method. Before starting it is necessary to give some preliminary definitions.
\subsection{Preliminary Definitions}
\begin{figure}
\centering
\begin{subfigure}[b]{0.2\textwidth}
%\captionsetup{justification=centering}
\caption{$x<p$} \label{fig:s01}
%\centering
\includegraphics[width=1.0\textwidth]{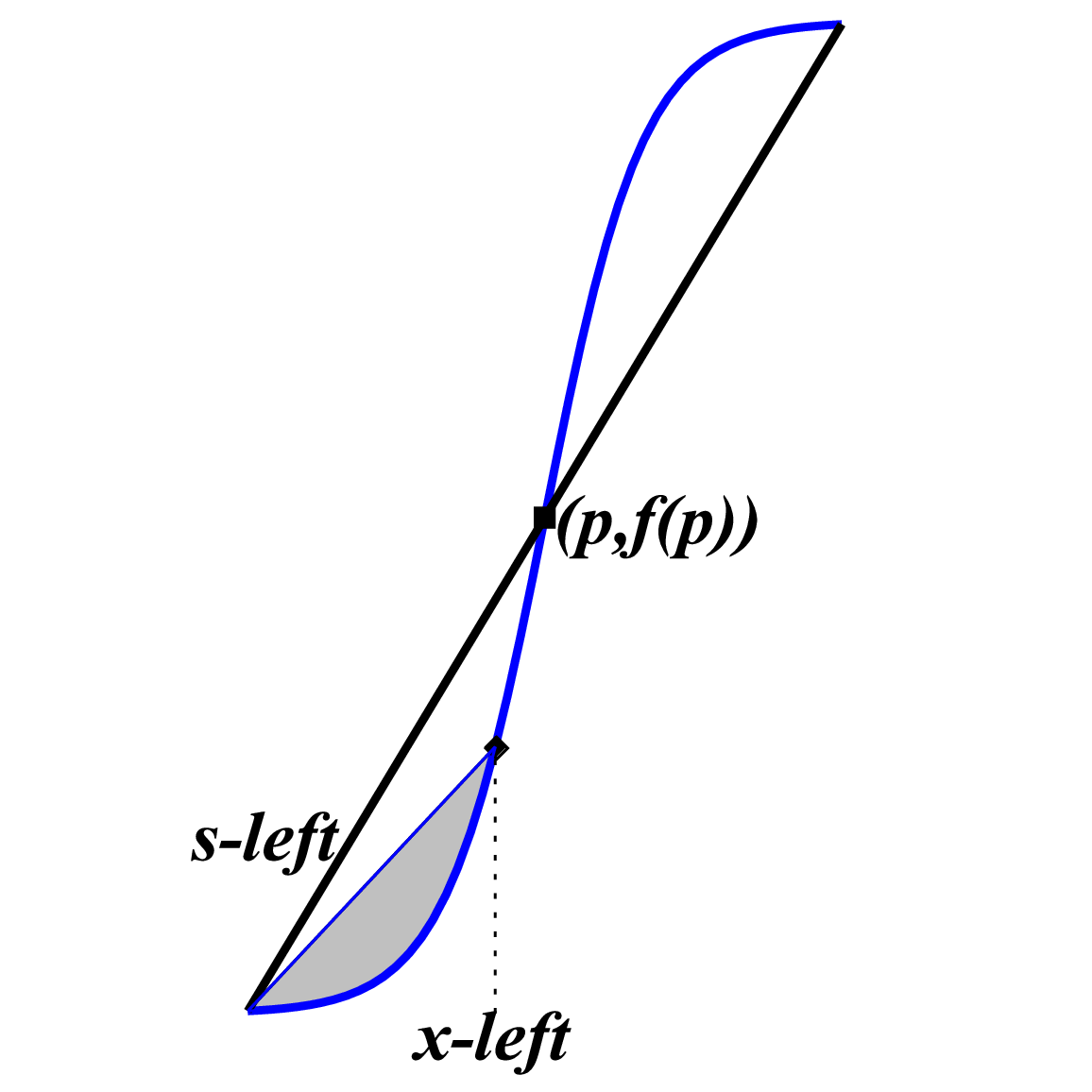}
\end{subfigure}
\begin{subfigure}[b]{0.2\textwidth}
%\captionsetup{justification=centering}
\caption{$x=p$} \label{fig:s02}
%\centering
\includegraphics[width=1.0\textwidth]{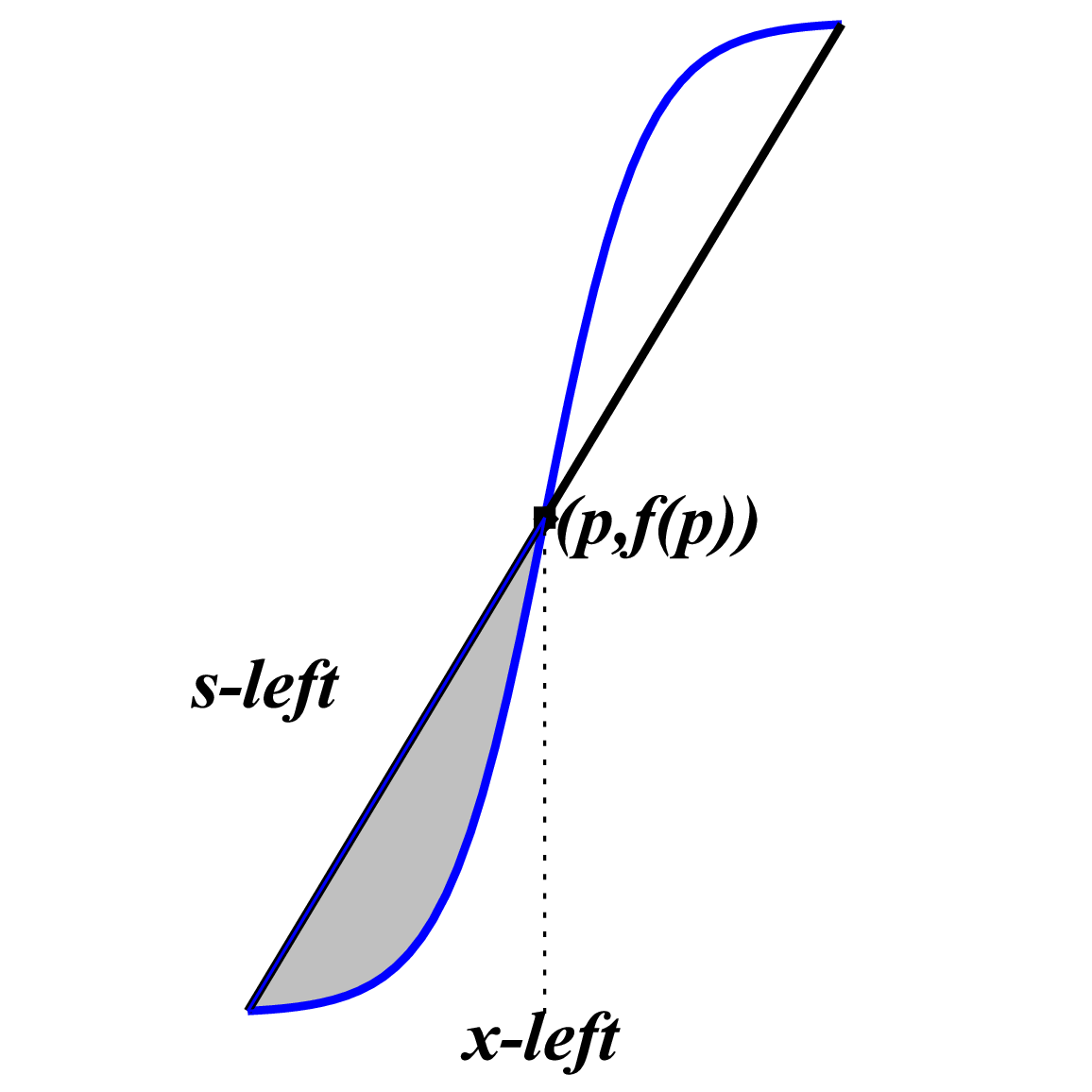}
\end{subfigure}
\begin{subfigure}[b]{0.2\textwidth}
%\captionsetup{justification=centering}
\caption{$x=x_{l}$} \label{fig:s03}
%\centering
\includegraphics[width=1.0\textwidth]{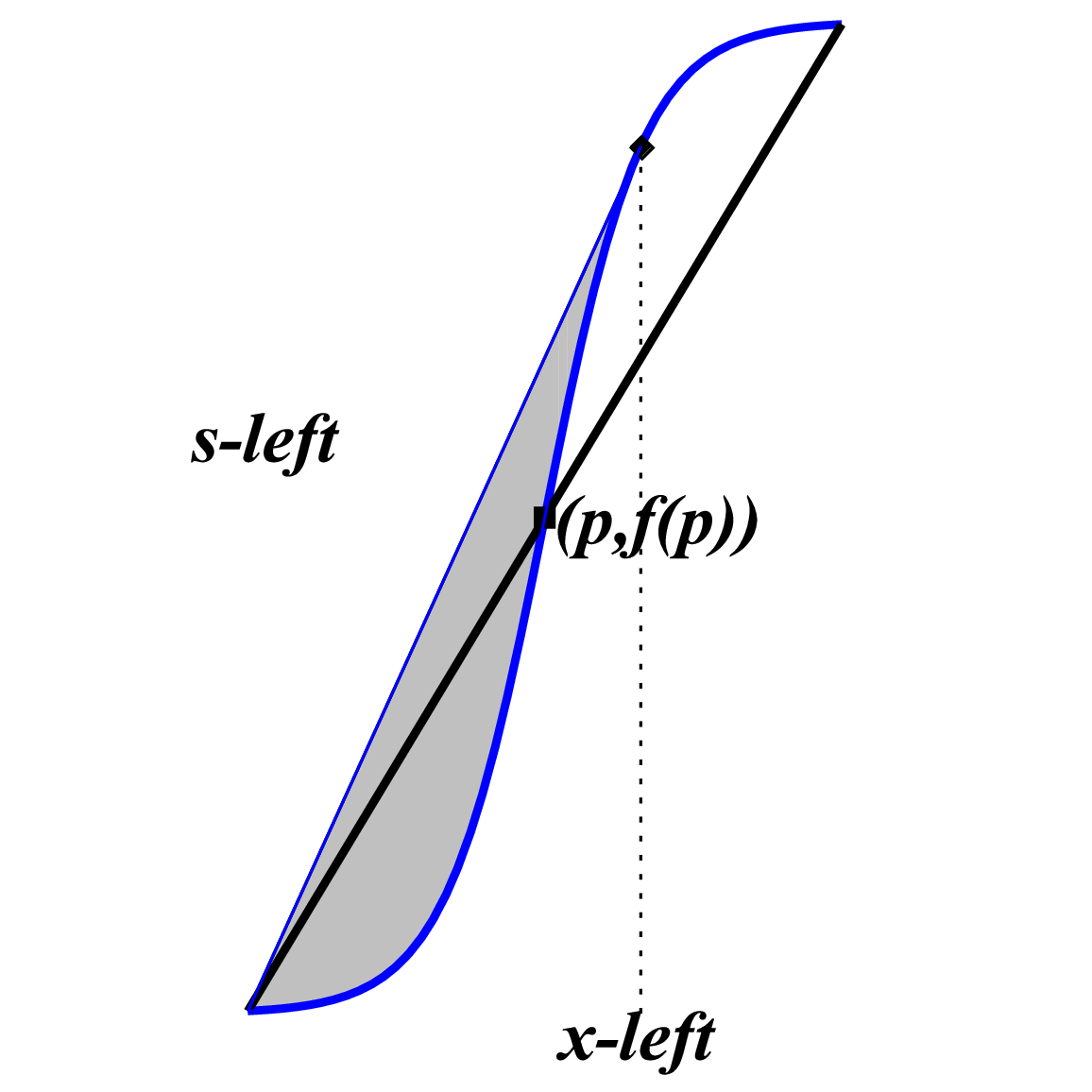}
\end{subfigure}
\begin{subfigure}[b]{0.2\textwidth}
%\captionsetup{justification=centering}
\caption{$x>x_{l}$} \label{fig:s04}
%\centering
\includegraphics[width=1.0\textwidth]{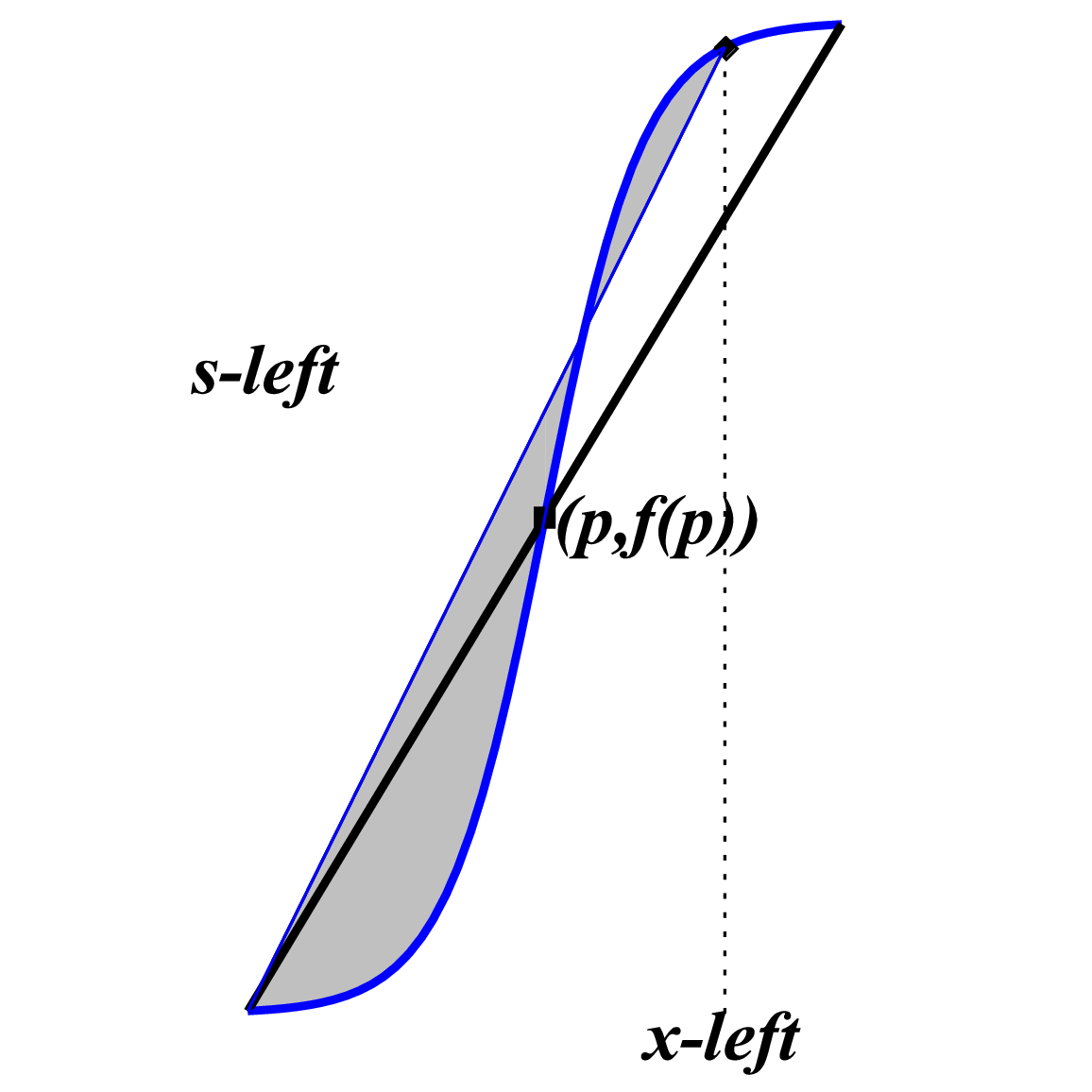}
\end{subfigure} \\ 
\begin{subfigure}[b]{0.2\textwidth}
%\captionsetup{justification=centering}
\caption{$x>p$} \label{fig:r01}
%\centering
\includegraphics[width=1.0\textwidth]{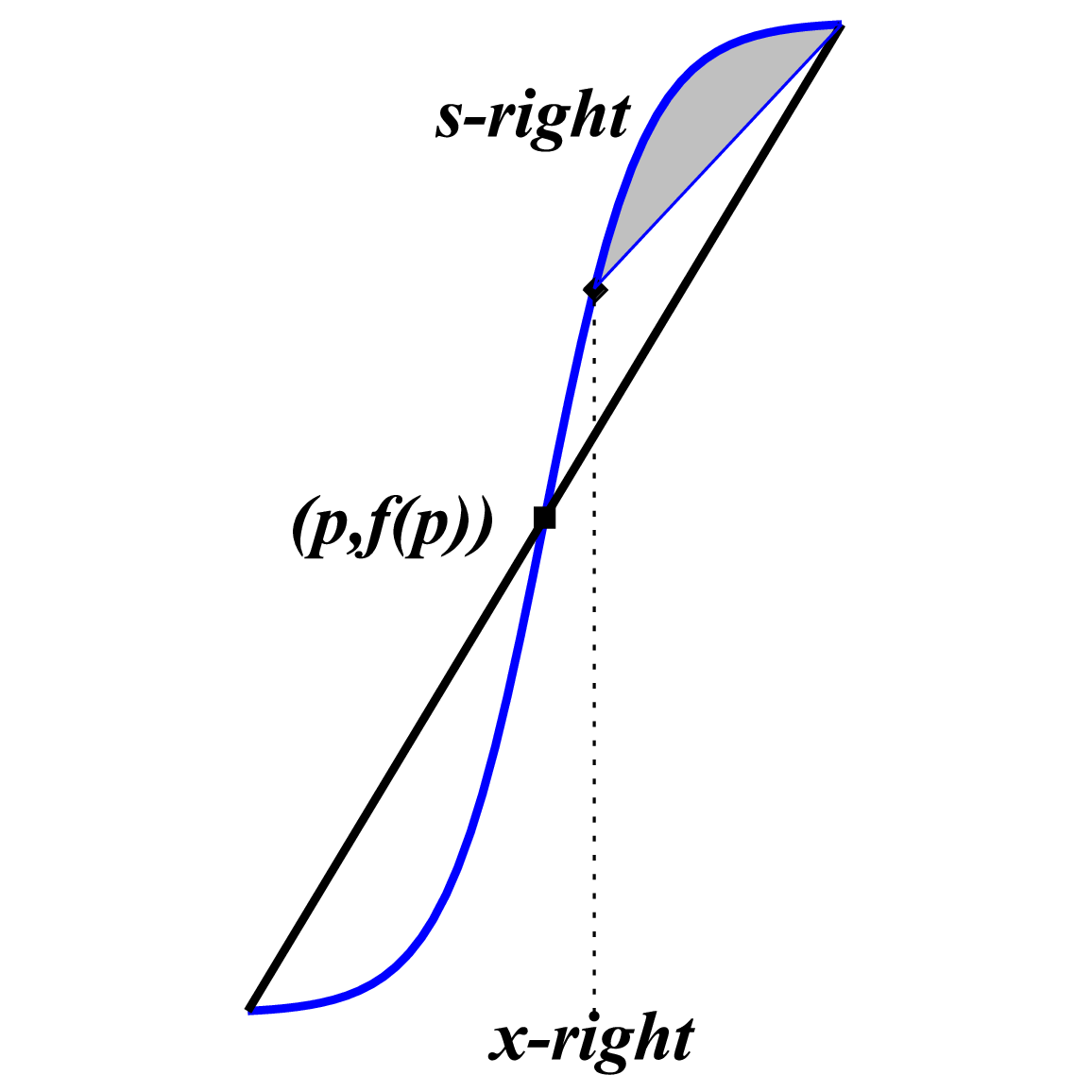}
\end{subfigure}
\begin{subfigure}[b]{0.2\textwidth}
%\captionsetup{justification=centering}
\caption{$x=p$} \label{fig:r02}
%\centering
\includegraphics[width=1.0\textwidth]{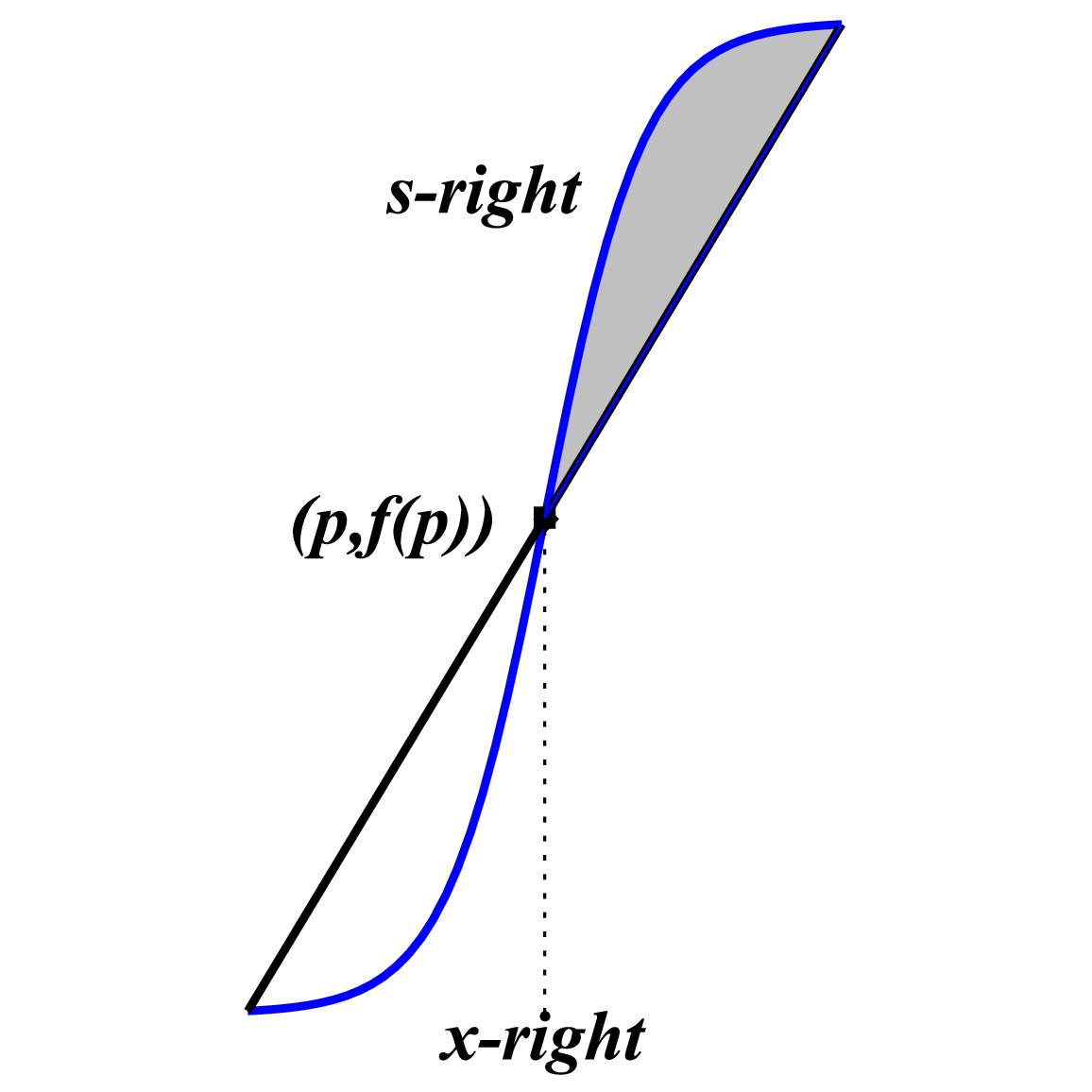}
\end{subfigure}
\begin{subfigure}[b]{0.2\textwidth}
%\captionsetup{justification=centering}
\caption{$x=x_{r}$} \label{fig:r03}
%\centering
\includegraphics[width=1.0\textwidth]{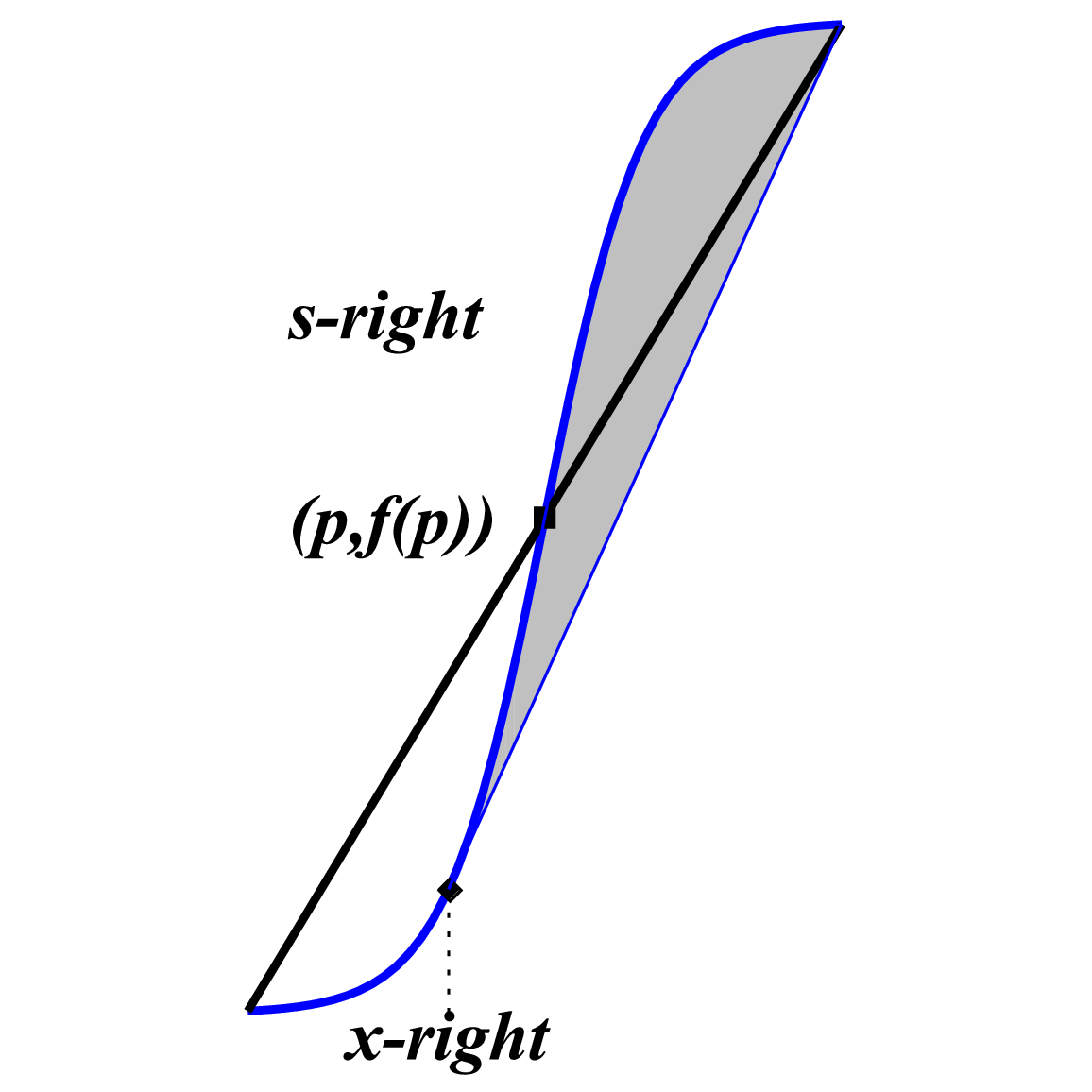}
\end{subfigure}
\begin{subfigure}[b]{0.2\textwidth}
%\captionsetup{justification=centering}
\caption{$x<x_{r}$} \label{fig:r04}
%\centering
\includegraphics[width=1.0\textwidth]{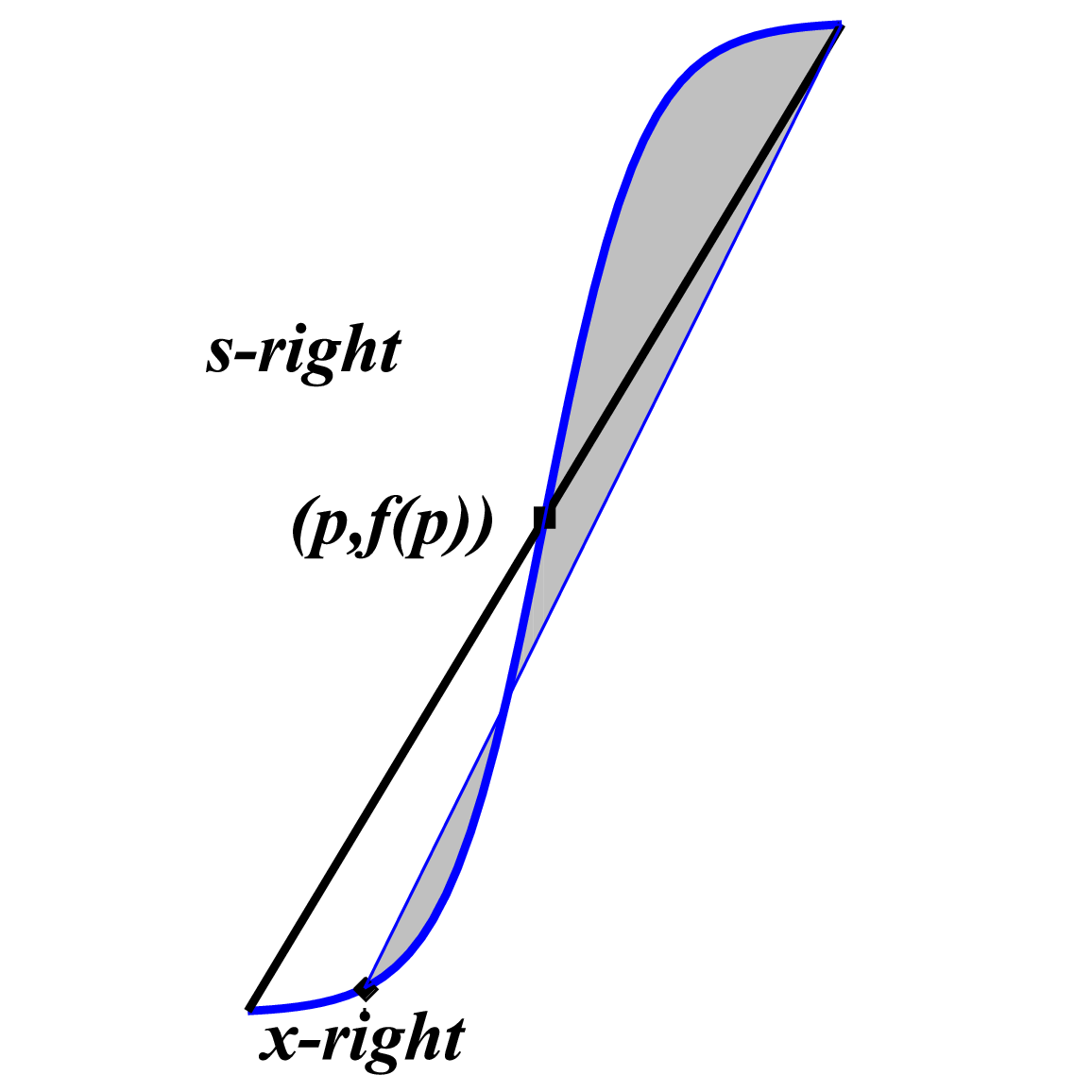}
\end{subfigure}
\caption{(color online) Illustration of the algebraic surfaces $s_{l}(a,x)<0$ (a-d), $s_{r}(b,x)>0$ (e-h) and of their corresponding points $x_{l},\,x_{r}$.}\label{fig:ese1}
\end{figure}

Let a function $f:[a,b]\rightarrow{R},\,\,f\in{C^{(n)}},\,n\ge{2}$ which is convex for $x\in[a,p]$ and concave for $x\in[p,b]$, p is the unique inflection point of $f$ in $[a,b]$ and let an arbitrary $x\in[a,b]$.
\begin{definition}\label{def:chords}
Total, left and right chord are the lines connecting points \\ 
$\{(a,f(a)), (b,f(b))\}$, $\{(a,f(a)), (x,f(x))\}$ and $\{(x,f(x)), (b,f(b))\}$ with Cartesian equations $g(x)$, $l(x)$ and $r(x)$ respectively.  
\end{definition}
If $cx+dy+e=0$ is the equation of the total chord, then by using elementary Analytical Geometry methods we can prove that the coefficients are
\begin{equation}\label{coef:gx}
\begin{matrix}
c=\frac{f(b)-f(a)}{b-a} & d=-1 & e={\frac {b\,f \left( a \right)-a\,f \left( b \right)}{b-a}} \\
\end{matrix}
\end{equation}
\begin{definition}\label{def:distFs}
Distance from total, left and right chord are the functions $F,F_{l},F_{r}:[a,b]\rightarrow{R}$ with:
\begin{eqnarray}
F(x)=f(x)-g(x) \\
F_{l}(x)=f(x)-l(x) \\
F_{r}(x)=f(x)-r(x)
\end{eqnarray}
\end{definition}
So, recalling \ref{coef:gx}, the distance from total chord is equal to
\begin{equation}
F \left( x \right) =f \left( x \right) -{\frac { \left( f \left( b
 \right) -f \left( a \right)  \right) x}{b-a}}-{\frac {b\,f \left( a \right)-a\,f \left( b \right)}{b-a}}
\end{equation}
\begin{definition}\label{def:slsr}
The s-left $(s_{l}(a,x))$ and s-right $(s_{r}(b,x))$ are the algebraic surfaces:
\begin{eqnarray}
s_{l}(a,x)=\int_{a}^{x}{F_{l}(t)dt} \\
s_{r}(b,x)=\int_{x}^{b}{F_{r}(t)dt} 
\end{eqnarray}
\end{definition}
\begin{definition}\label{def:xlxr}
The x-left $(x_{l})$ and x-right $(x_{r})$ are values such that:  
\begin{eqnarray}
x_{l}=\underset{x\in{[a,b+\delta_1]}}{argmin}\{s_{l}\left(a,x\right)\} \\
x_{r}=\underset{x\in{[a-\delta_2,b]}}{argmax}\{s_{r}\left(b,x\right)\} 
\end{eqnarray}
with $\delta_{1},\delta_{2}>0$ taken as small as necessary for $x_{l},x_{r}$ to be unique unconstrained extremes in the corresponding intervals. 
\end{definition}
A graphical illustration of the above defined left and right-terms is presented at Figure \ref{fig:ese1} where we observe that when $x=x_{l},\,x=x_{r}$ then we achieve the algebraic minimum $s_{l}(a,x)$ and maximum $s_{r}(b,x)$, respectively. The motivation for the \emph{left-} and \emph{right-} naming was due to the origin of the chords: left chord starts from the left edge while right chord starts from the right edge of the graph. \\
%\end{definition}
\begin{definition}\label{def:normdist}
Normal algebraic distance of the curve point $(x,f(x))$ from the total chord \ref{def:chords} \& \ref{coef:gx} is the function $N:[a,b]\rightarrow{R}$ with:
\begin{equation}
N(x)=-\frac{cx+df(x)+e}{\sqrt{c^2+d^2}} 
\end{equation}
\end{definition}
For a convex/concave curve the above definition gives $N(x)<0$ when $x<p$ and $N(x)>0$ if $x>p$. Additionally, for computation purposes, by using \ref{coef:gx}, the above distance is 
\begin{equation}
N \left( x \right) =-{\frac {  \left( f \left( b \right) -f
 \left( a \right)  \right) x- \left( b-a \right) f \left( x \right) +b\,f
  \left( a \right)-a\,f \left( b \right)   }{\sqrt { \left( f
 \left( b \right) -f \left( a \right)  \right) ^{2}+ \left( b-a
 \right) ^{2}}}}
\end{equation}

We call \emph{standard partition} (\emph{SP}) the strictly sorted grid of points, not necessary equal spaced:
\beq\label{spxi}
\left\{x_{i},i=0,1,\ldots{n},\,\,a=x_0<x_1<\ldots<x_n=b\right\}
\eeq 
The corresponding $(x_i,f_i)$ data produced from $f$ by the no error process:
\beq
\label{eq:noerrdat}
\left\{f_{i}=f(x_{i}),i=0,1,\ldots{n}\right\}
\eeq
If errors occur then we have the $(x_i,\phi_i)$ noisy data produced from $f$ by the process:
\beq
\label{eq:errdat}
\left\{\phi_{i}=f(x_{i})+\epsilon_{i},\,\,\epsilon_{i}\sim{iid(0,\sigma^2)}\,\,i=0,1,\ldots{n}\right\}
\eeq
Our analysis is focused on uniform distributions $(\epsilon_{i}\sim{U(-r,r)})$ but it is applicable for every distribution with zero mean, for example the normal $N(0,\sigma^2)$. If the error distribution is not a zero mean one, then the results are ambiguous.\\
\begin{definition}\label{Phis}
For the noisy data \ref{eq:errdat} we define the discrete distances from total, left and right chord as the values
\begin{eqnarray}
\Phi(x_i)=\phi(x_i)-g(x_i),\,\,i=0,1,\ldots,n \\
\Phi_{l}(x_i)=\phi(x_i)-l(x_i),\,\,i=0,1,\ldots,n \\
\Phi_{r}(x_i)=\phi(x_i)-r(x_i),\,\,i=0,1,\ldots,n 
\end{eqnarray}
\end{definition}

\begin{definition}\label{def:symp}
Function $f$ is called symmetric around inflection point or \emph{symmetry around inflection point} exists when:
\beq 
f \left( p+x \right) -f \left( p \right) =f\left( p \right)-f\left( p-x \right),\forall{x}\in{R}
\eeq 
or
\beq
f \left( p+x \right)+f\left( p-x \right)-2\,f\left( p \right)=0,\,\forall{x}\in{R}
\eeq
\end{definition}
\begin{definition}\label{def:symped}
Function $f$ is called locally $\left(\epsilon,\delta\right)$ asymptotically symmetric around inflection point or 
\emph{local $\epsilon-\delta$ asymptotic symmetry} exists when:
\beq
\left|f(p+x)+f(p-x)-2f(p)\right|<\epsilon,\,\forall{x}\in(p-\delta,p+\delta)
\eeq
\end{definition}

\begin{definition}\label{def:asymlr}
For a function $f:[a,b]\rightarrow{R}$ we have that:
\begin{enumerate}
  \item Data symmetry w.r.t. inflection point exists if $p-b=a-p$
	\item Data left asymmetry w.r.t. inflection point exists if $p-b<a-p$
	\item Data right asymmetry w.r.t. inflection point exists if $p-b>a-p$	
\end{enumerate}
\end{definition}

\begin{definition}\label{def:totsym}
A function $f:[a,b]\rightarrow{R}$ is called totally symmetric or total symmetry exist, 
if it is symmetric around inflection point p and also exist data symmetry w.r.t. p.
\end{definition}

\begin{definition}\label{def:trapezrule}
For every subsequent $x_i<x_j$ the elementary trapezoidal estimation holds: 
\beq
\int_{x_i}^{x_j}{f(x)dx} \approx T_{i,j}(f,x_i,x_j)=\frac{f(x_i)+f(x_j)}{2}(x_j-x_i)
\eeq
And for every standard partition the total trapezoidal estimation holds:
\begin{equation}
\int_{a}^{b}{f(x)dx} \approx T_{n+1}(f,a,b)=\sum_{i=0}^{n-1} T_{i,i+1}(f,x_i,x_{i+1})  
\end{equation}
\end{definition}

\subsection{The Extremum Surface Method}  
We can prove that:
\begin{lemma}\label{lem:xlxr}
The x-left $(x_{l})$ and x-right $(x_{r})$ are the points where left and right chord respectively are tangent to the graph $G(f)$.
\end{lemma}
\textbf{Proof}\\
\begin{enumerate}
	\item Let $x_l\leq{p}$ be the first point where $l(x)$ cuts $G(f)$ from left to right. The function is convex for $x\in[a,p]$, so $G(f)$ is always below the left chord, 
	thus giving a negative value for $s_{l}\left(a,x_l\right)$, 
	which is increasing in absolute value as $x_l$ departures from a to the right until point p. After inflection point f is increasing, so it is possible to continue 
	adding negative values of surface until the point $x^{*}\in(p,b]$ for which $G(f)$ and $l(x)$ have one only common point.\\
	If we continue beyond this point, then $G(f)$ and $l(x)$ have again two common points	$(x_1,y_1)$, \\
	$(x_2,y_2)$, $x_1<x_2$ such that $s_{l}\left(a,x_1\right)<0$ and $s_{l}\left(x_1,x_2\right)>0$, thus we have started adding positive values to the $s_{l}\left(a,x\right)$ and this leads to 
	a raise of the total value $s_{l}\left(a,x\right)=s_{l}\left(a,x_1\right)+s_{l}\left(x_1,x_2\right)$.\\
	So, function $s_{l}\left(a,x\right)$ is decreasing for $x\in[a,x^{*}]$ and increasing for $x\in[x^{*},b]$, thus $x^{*}>p$ is a local minimum and we call it $x_l$. 
  \item Let $x_r\ge{p}$ be the first point where $r(x)$ cuts $G(f)$ from right to left. The function is concave for $x\in[p,b]$, so $G(f)$ is always above the right chord, 
  thus giving a positive value for $s_{r}\left(b,x_r\right)$,
  which is increasing as $x_r$ departures from b to the left until inflection point p. After that point f is still above the right chord until the point $x^{*}\in[a,p)$ for 
  which $G(f)$ and $r(x)$ have only one common point.\\
  If we continue again beyond this point, then $G(f)$ and $r(x)$ have again two common points $(x_1,y_1),(x_2,y_2),\,\,x_1>x_2$ such that 
  $s_{r}\left(b,x_1\right)>0$ and $s_{r}\left(x_1,x_2\right)<0$, so we have started adding negative values to the $s_{r}\left(b,x\right)$ and this leads to 
  a reduction of the total value $s_{r}\left(b,x\right)=s_{r}\left(b,x_1\right)+s_{r}\left(x_1,x_2\right)$.\\
  So, function $s_{r}\left(b,x\right)$ is increasing for $x\in[a,x^{*}]$ and decreasing for $x\in[x^{*},b]$, thus $x^{*}<p$ is a local minimum and we call it $x_r$.  
\end{enumerate}
%\qquad
\hfill\(\Box\) \\

\begin{corollary}\label{cor:xtanlr}
For the definitions of \ref{def:xlxr} it holds
\begin{eqnarray}
x_{l}=\underset{x\in{[a,b+\delta_1]}}{arg}\left\{f^{'}\left(x\right)=\frac{f(x)-f(a)}{x-a}\right\} \\
x_{r}=\underset{x\in{[a-\delta_2,b]}}{arg}\left\{f^{'}\left(x\right)=\frac{f(b)-f(x)}{b-x}\right\}
\end{eqnarray}
with $\delta_{1},\delta_{2}>0$ taken as small as necessary for $x_{l},x_{r}$ to be unique unconstrained solutions in the corresponding intervals. 
\end{corollary}
This tangency condition is obvious from Figure \ref{fig:s03}.\\
\begin{corollary}\label{lem:axrxlb}
Let a function $f:[a,b]\rightarrow{R},\,\,f\in{C^{(n)}},\,n\ge{2}$ which is convex for $x\in[a,p]$ and concave for $x\in[p,b]$. Then we have one of the following possibilities:
\begin{enumerate}
	\item If $x_{l},x_{r}\in{[a,b]}$ then $a\leq{x_{r}}<x_{l}\leq{b}$
	\item If $x_{l}\notin{[a,b]}$ then $x_{l}>b$
	\item If $x_{r}\notin{[a,b]}$ then $x_{r}<a$
\end{enumerate}
\end{corollary}
We define the next theoretical estimator of the inflection point:
\begin{definition}\label{def:tese}
The theoretical extremum surface estimator (TESE) is
\beq
x_{S}=
 \begin{Bmatrix}
 \frac{x_{l}+x_{r}}{2} &{,}&x_{l},x_{r}\in{[a,b]} \\
 \frac{b+x_{r}}{2} &{,}&x_{l}>b \\
 \frac{x_{l}+a}{2} &{,}&x_{r}<a 
 \end{Bmatrix}   
\eeq
\end{definition}

\begin{lemma}\label{lem:trapezest}
If the mesh $\lambda(n)$ of the standard partition is such that:
\begin{displaymath}
\lim_{n \to \infty}\,n\lambda(n)^2=0
\end{displaymath}
then $T_{n+1}(\phi,a,b)$ is a consistent estimator of the true value of $T_{n+1}(f,a,b)$.
\end{lemma}
\textbf{Proof}\\
For every subsequent $x_i<x_{i+1}$ the elementary trapezoidal estimation is: 
\beq
T_{i,i+1}(\phi,x_i,x_{i+1})=\frac{x_{i+1}-x_i}{2}\phi(x_i)+\frac{x_{i+1}-x_i}{2}\phi(x_{i+1})
\eeq
Taking the expected value we obtain:
\beq
\begin{array}{lll}
E\left(T_{i,i+1}(\phi,x_i,x_{i+1})\right)&=&\frac{x_{i+1}-x_i}{2}E\left(\phi(x_i)\right)+\frac{x_{i+1}-x_i}{2}E\left(\phi(x_{i+1})\right) \\
&=&\frac{x_{i+1}-x_i}{2}f(x_i)+\frac{x_{i+1}-x_i}{2}f(x_{i+1}) \\
&=&T_{i,i+1}(f,x_i,x_{i+1}) 
\end{array}
\eeq
so from the linearity of expected value we have also that:
\beq
E\left(T_{n+1}(\phi,a,b)\right)=\sum_{i=0}^{n-1} E\left(T_{i,i+1}(\phi,x_i,x_{i+1})\right)=T_{n+1}(f,a,b)
\eeq
Thus our estimator is unbiased. \\
We continue by computing the variance of the elementary trapezoidal estimation:
\beq
\begin{array}{lll}
V\left(T_{i,i+1}(\phi,x_i,x_{i+1})\right)& =& \left(\frac{x_{i+1}-x_i}{2}\right)^{2}V\left(\phi(x_i)\right)+\left(\frac{x_{i+1}-x_i}{2}\right)^{2}V\left(\phi(x_{i+1})\right) \\ 
& =& \frac{(x_{i+1}-x_i)^2}{4}\sigma^2+\frac{(x_{i+1}-x_i)^2}{4}\sigma^2 \\ 
& =&\frac{(x_{i+1}-x_i)^2}{2}\sigma^2
\end{array}
\eeq
We have two cases.\\
If standard partition is equal spaced, then $x_{i+1}-x_i=\frac{b-a}{n}$ and we obtain:
\beq
V\left(T_{i,i+1}(\phi,x_i,x_{i+1})\right)=\frac{(b-a)^2}{2{n^2}}\sigma^2
\eeq
Let' s compute now the variance of estimator $T_{n+1}(\phi,a,b)$:
\beq
\begin{array}{lll}
V\left(T_{n+1}(\phi,a,b)\right)&=&V\left(\sum\limits_{i=0}^{n-1} T_{i,i+1}(\phi,x_i,x_{i+1})  \right) \\
&=&n\,V\left(T_{i,i+1}(\phi,x_i,x_{i+1})\right) \\
&=&n\,\frac{(b-a)^2}{2{n^2}}\sigma^2 \\
&=&\frac{(b-a)^2}{2n}\sigma^2 \\
\end{array}
\eeq	
so clearly it holds:
\begin{displaymath}
\lim_{n \to \infty}\,V\left(T_{n+1}(\phi,a,b)\right)=\lim_{n \to \infty}\frac{(b-a)^2}{2n}\sigma^2=0
\end{displaymath}
For the second case, if standard partition is not equal spaced then the mesh or norm of the partition is 
\begin{displaymath}
\lambda(n)=\max_{i=0,\ldots,{n-1}}(x_{i+1}-x_i)
\end{displaymath}
Then it is easy to show that:
\beq
V\left(T_{i,i+1}(\phi,x_i,x_{i+1})\right)\leq\frac{\lambda(n)^2}{2}\sigma^2
\eeq
and the total variance is:
\beq
V\left(T_{n+1}(\phi,a,b)\right)\leq \frac{\sigma^2}{2}\,n\lambda(n)^2 \underset{n \to \infty}{\longrightarrow} {0}
\eeq
from our hypothesis.So the estimator is consistent. \hfill\(\Box\) \\
Now we are able to compute using our trapezoidal rule \ref{def:trapezrule} data estimations for $s_l(x_0,x_j)$ and $s_r(x_n,x_j)$:
\begin{definition}\label{def:slest}
The data estimators of the algebraic surfaces of \ref{def:slsr} are
\begin{eqnarray}
s_{l,j+1}(x_0,x_j)=T_{j+1}(\Phi_{l},x_{0},x_{j}) \\
s_{r,n-j+1}(x_n,x_j)=T_{n-j+1}(\Phi_{r},x_{j},x_{n}) \\
\end{eqnarray}
\end{definition}
It is time to define the next data estimators for $x_l,x_r$.
\begin{definition}\label{def:chilrer}
The $\chi_{l},\,\,\chi_{r}$ are values such that:  
\begin{eqnarray}
\chi_{l}=x_{j_l}  \\ \nonumber
j_{l}=\underset{j\in{[1,n]}}{argmin}\{s_{l,j+1}\left(x_0,x_j\right)\} \\ 
%j_l=argmin_{j\in{[1,n]}}\{s_{l,j+1}\left(x_0,x_j\right)\} \\ 
\chi_{r}=x_{j_r} \\ \nonumber
j_{r}=\underset{j\in{[0,n-1]}}{argmax}\{s_{r,n-j+1}\left(x_n,x_j\right)\} \\
%j_r=argmax_{j\in{[0,n-1]}}\{s_{r,n-j+1}\left(x_n,x_j\right)\} \\
\end{eqnarray}
\end{definition}
We define now the noisy data estimator of the inflection point:
\begin{definition}\label{def:ese}
The data extremum surface estimator (ESE) is
\beq
\chi_{S}=\frac{\chi_{l}+\chi_{r}}{2}
\eeq
\end{definition}

\begin{lemma}\label{lem:esecon}
The ESE is a consistent estimator of TESE with all relevant integrals 
calculated via trapezoidal rule.
\end{lemma}
\textbf{Proof}\\
We have proven in Lemma \ref{lem:trapezest} that trapezoidal rule for the noisy data gives a consistent estimator for 
the trapezoidal estimation of the actual data, thus $\chi_l,\chi_r$ are consistent estimators of the true  $x_l,x_r$ respectively, with relevant integrals trapezoidal calculated.\\
If the interval $[a,b]$ is such that both $x_l,x_r\in[a,b]$ then ESE is a consistent estimator of trapezoidal calculated $x_{S}=\frac{x_l+x_r}{2}$.\\
If $x_l>b$ then recalling Proof of Lemma \ref{lem:xlxr} $s_{l}\left(a,x\right)$ is a decreasing function, so the minimum $\chi_l$ is achieved when $\chi_l=b$, the rightmost value of $[a,b]$.\\  
If $x_r<a$ then recalling Proof of Lemma \ref{lem:xlxr} $s_{r}\left(b,x\right)$ is an increasing function, so the maximum $\chi_r$ is achieved when $\chi_r=a$, the leftmost value of $[a,b]$.\\ 
So, for every possible case, ESE is a consistent estimator of the TESE given by integrals calculated via trapezoidal rule. %\qquad 
 \hfill\(\Box\) \\ 
We have to make a remark about the concept of convex or concave area, as is defined in \cite{mok-08} and as defined in this work. There the concept of area that is computed is a summation of the distances from a chord and curve, between two critical points, while our approach computes an actual geometrical area by using the trapezoidal rule.

\subsection{The Extremum Distance Method}  
\begin{figure}[hbtp]
\begin{center}
\includegraphics[scale=0.3,angle=0]{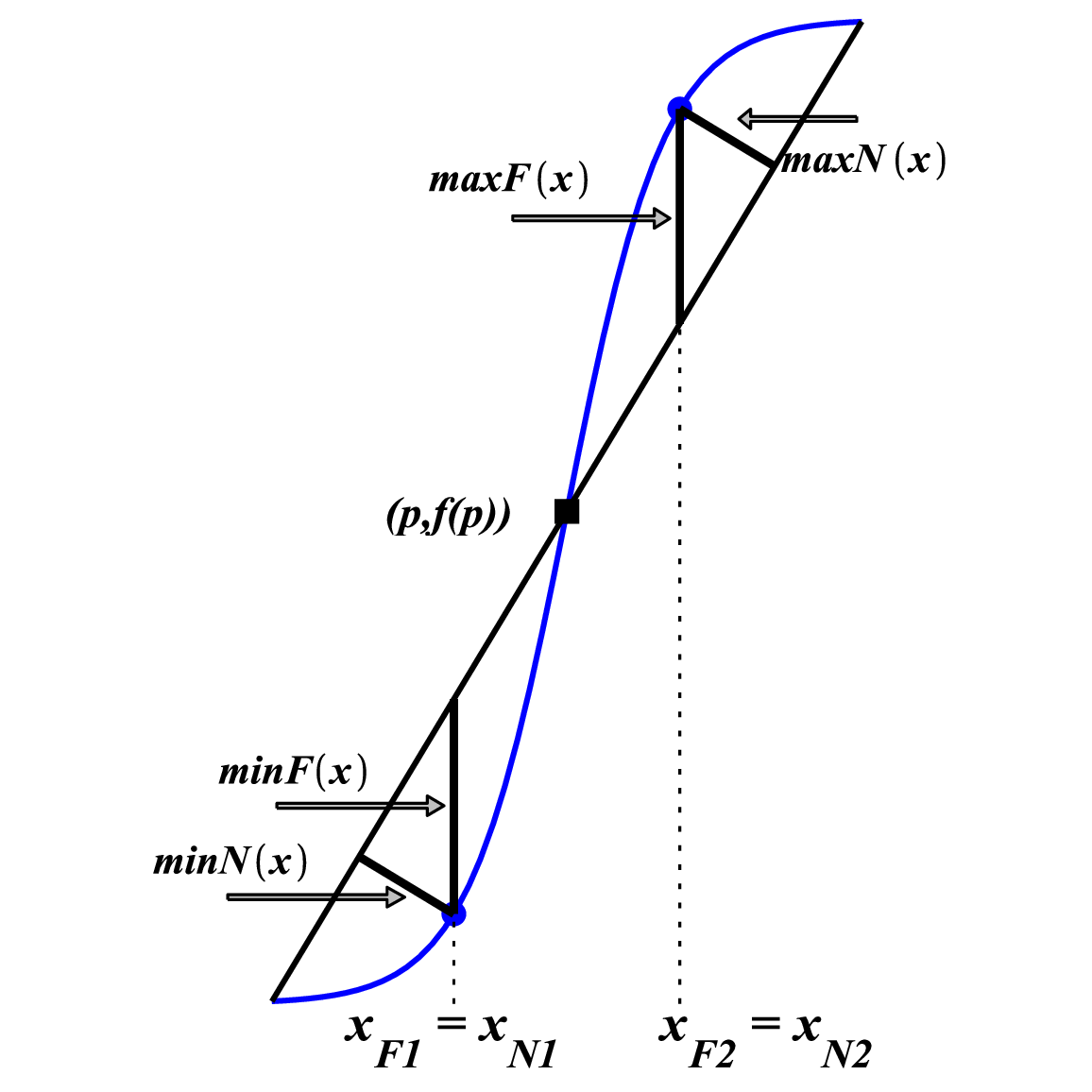}
\caption{(color online) Illustration of points $x_{F1},\,x_{F2}$} and the corresponding distances ${min}{F(x)},\,{max}{F(x)}$ \label{fig:pede1}
\end{center}
\end{figure}

\begin{definition}\label{def:xf12}
The xF-left $(x_{F1})$, xF-right $(x_{F2})$ and xN-left $(x_{N1})$, xN-right $(x_{N2})$ are values such that:
\begin{eqnarray}
x_{F1}=\underset{x\in{[a-\delta_1,b]}}{argmin}{F(x)},\quad x_{F2}=\underset{x\in{[a,b+\delta_2]}}{argmax}{F(x)} \\
x_{N1}=\underset{x\in{[a-\delta_1,b]}}{argmin}{N(x)},\quad
x_{N2}=\underset{x\in{[a,b+\delta_2]}}{argmax}{N(x)} 
\end{eqnarray}
with $\delta_{1},\delta_{2}>0$ taken as small as necessary for $x_{F1},x_{F2}$ to be unique unconstrained extremums in the corresponding intervals. 
\end{definition}
The defined points and the corresponding line segments can be found at the Figure \ref{fig:pede1} where it is also obviously shown that when we achieve a stationary value for $F(x)$ of \ref{def:distFs}, then we achieve also the relevant stationary value for normal distance $N(x)$ of \ref{def:normdist}, since the vector defined from $N(x)$ is just the orthogonal projection of the vector defined from $F(x)$ at the normal vector to the total chord. Now we shall prove the next useful Lemma.\\
 
\begin{lemma}\label{lem:xf12}
For the definitions of \ref{def:xf12} it holds
\beq
{x_{F}}_{1,2}=\underset{x\in[a-\delta_1,b+\delta_2]}{arg}\left\{f^{'}\left(x\right)=\frac{f(b)-f(a)}{b-a}\right\}
%{x_{F}}_{1,2}=arg_{x\in[a-\delta_1,b+\delta_2]}\{f^{'}\left(x\right)=\frac{f(b)-f(a)}{b-a}\}
\eeq
with $\delta_{1},\delta_{2}>0$ taken as small as necessary for $x_{F1},x_{F2}$ to be unique unconstrained extremes in the corresponding intervals.
\end{lemma}
\textbf{Proof}\\
We have extended the interval $[a,b]$ such that there exist (both unconstrained) a local minimum and a local maximum inside. For our convex/concave case let $\rho\in[a-\delta_1,b+\delta_2],\rho\notin\{a,b\}$ is the internal root of $F(x)$. Then we have that $F(x)<0,x\in[a-\delta_1,\rho]$ and $F(x)>0,x\in[\rho,b+\delta_2]$, because function is convex near a and concave near b. Thus the local minimum exists at $[a-\delta_1,\rho]$ and the local maximum lies in $[\rho,b+\delta_2]$. If we take the first derivative we have that:
\beq
F^{ '}(x)=f^{ '}(x)-g^{ '}(x)=f^{ '}(x)-\lambda=0\Rightarrow f^{ '}(x)=\lambda
\eeq
where $\lambda=\frac{f(b)-f(a)}{b-a}$ is the slope of the total chord. But the above equation must hold for both local minimum/maximum ${x_{F}}_{1,2}$, so it is necessary to hold:
\beq
f^{ '}({x_{F}}_{1})=f^{ '}({x_{F}}_{2})=\frac{f(b)-f(a)}{b-a}
\eeq 
We can also check the second derivative which is:
\beq
F^{ ''}(x)=f^{ ''}(x)
\eeq
so it holds $F^{ ''}({x_{F}}_{1})=f^{ ''}({x_{F}}_{1})>0$ and $F^{ ''}({x_{F}}_{2})=f^{ ''}({x_{F}}_{2})<0$, i.e. we have the correct signs for local minimum and maximum respectively.
\hfill\(\Box\) \\
\begin{corollary}\label{cor:axfb}
Let a function $f:[a,b]\rightarrow{R},\,\,f\in{C^{(n)}},\,n\ge{2}$ which is convex for $x\in[a,p]$ and concave for $x\in[p,b]$. Then we have one of the following possibilities:
\begin{enumerate}
	\item If $x_{F1},x_{F2}\in{[a,b]}$ then $a\leq{x_{F1}}<{x_{F2}}\leq{b}$
	\item If $x_{F1}\notin{[a,b]}$ then $x_{F1}<a$
	\item If $x_{F2}\notin{[a,b]}$ then $x_{F2}>b$
\end{enumerate}
\end{corollary}
We define the next theoretical estimator of the inflection point:
\begin{definition}\label{def:tede}
The theoretical extremum distance from total chord estimator (TEDE) is such that
\beq
x_{D}=\frac{x_{F1}+x_{F2}}{2}
\eeq
\end{definition}

Now we can define the data estimators of $x_{F1},x_{F2}$ and $x_{D}$.
\begin{definition}\label{def:chiF12}
The data estimations of the points defined at \ref{def:xf12} are 
\begin{eqnarray}
\chi_{F1}=x_{j_1} \\ \nonumber
j_1=\underset{j\in{[0,n]}}{argmin}\{\Phi(x_j)\} \\
%j_1=argmin_{j\in{[0,n]}}{\Phi(x_j)} \\
\chi_{F2}=x_{j_2} \\ \nonumber
j_2=\underset{j\in{[0,n]}}{argmax}\{\Phi(x_j)\}
%j_2=argmax_{j\in{[0,n]}}{\Phi(x_j)}
\end{eqnarray}
\end{definition}
\begin{definition}\label{def:ede}
The extremum distance from total chord estimator (EDE) is
\beq
\chi_{D}=\frac{\chi_{F1}+\chi_{F2}}{2}\,\,\text{iff}\,\,\chi_{F2}\geq\chi_{F1}
\eeq
\end{definition}
\begin{lemma}\label{lem:unbiasedede}
The EDE is an unbiased estimator of TEDE.
\end{lemma}
\textbf{Proof}\\
For all $\Phi(x_j),j=0,1,\ldots,n$ it holds:
\beq
E\left(\Phi(x_j)\right)=F(x_j)
\eeq
so if we take the noisy data instead of the true data it has to be also that:
\begin{eqnarray}
E\left(\underset{j\in{[0,n]}}{min}\{\Phi(x_j)\}\right)=\underset{j\in{[0,n]}}{min}\{F(x_j)\} \\ \nonumber
E\left(\underset{j\in{[0,n]}}{max}\{\Phi(x_j)\}\right)=\underset{j\in{[0,n]}}{max}\{F(x_j)\}
\end{eqnarray}
\hfill\(\Box\) \\
At this stage we have to mention that there exists a similar work with distances from the chord, see \cite{han-01}, where a summation of the relevant distances from the chord is taken in order to define the concept of a \emph{discrete curvature} for a planar curve.\\
Similar approach is that of \cite{mok-08}, where it is also used a proper summation of the distances between the chord and the curve points. 
%Another method, see \cite{mok-86}, is using a Gaussian kernel in order to smooth first the curve and then to compute the curvature.\\
 Here we do not define and we do not compute any kind of curvature, but we just choose only the two extreme distances needed for Definition \ref{def:xf12}.

\subsection{Iterative application of geometrical methods}
\par Another very important opportunity is the possibility of iterations like the well known \emph{bisection method} in root finding. Recall that for a continuous function if $f(\alpha)\,f(\beta)<0$ then exist $\xi\in(\alpha,\beta)$ such that $f(\xi)=0$. Our ESE method always gives an interval that contains the true inflection point p or a point close to the edge a or b, if data is just convex (or just concave) and inflection point does not exist. EDE method also gives an interval in most cases, although it is more sensitive to errors, so it does not always give a point close to a or b, if simple convexity or concavity exist.
\begin{enumerate}
	\item 
\emph{ESE iterative method or Bisection-ESE or BESE}
\par We apply to initial data $\left\{(x_i,\phi_i),i=0,\ldots,n\right\}$ the ESE method and have the $0^{th}$ output for it:\\
$$
\begin{array}{l}
[j_r^{(0)},j_l^{(0)}],\,\chi_r^{(0)}=x_{j_r^{(0)}},\,\chi_l^{(0)}=x_{j_l^{(0)}},\,\chi_{S}^{(0)}=\frac{\chi_r^{(0)}+\chi_l^{(0)}}{2} \\
\end{array}
$$
If and only if $j_l^{(0)}>j_r^{(0)}$, then we apply again ESE for data:\\
$$
\left\{(x_i,\phi_i),i=j_r^{(0)},\ldots,j_l^{(0)}\right\}
$$
and obtain the $1^{st}$ output for it:\\
$$
\begin{array}{l}
[j_r^{(1)},j_l^{(1)}],\,\chi_r^{(1)}=x_{j_r^{(1)}},\,\chi_l^{(1)}=x_{j_l^{(1)}},\,\chi_{S}^{(1)}=\frac{\chi_r^{(1)}+\chi_l^{(1)}}{2} \\
\end{array}
$$
\par We continue until $j_l^{(k)}<j_r^{(k)}$ or until $\left|\chi_{S}^{(k)}-\chi_{S}^{(k-1)}\right|<e$, with $e=10^{-8}$ to be a good tolerance for all examined data. 
\item 
\emph{EDE iterative method or Bisection-EDE or BEDE}
\par We apply to initial data $\left\{(x_i,\phi_i),i=0,\ldots,n\right\}$ EDE method and have the $0^{th}$ output iff $x_{F2}>x_{F1}$:\\
$$
[j_1^{(0)},j_2^{(0)}],\,\chi_{F1}^{(0)}=x_{j_1^{(0)}},\,\chi_{F2}^{(0)}=x_{j_2^{(0)}},\,\chi_{D}^{(0)}=\frac{\chi_{F1}^{(0)}+\chi_{F2}^{(0)}}{2}
$$
If and only if $j_2^{(0)}>j_1^{(0)}$, then we apply again EDE method for data:\\
$$
\left\{(x_i,\phi_i),i=j_1^{(0)},\ldots,j_2^{(0)}\right\}
$$
and obtain the $1^{st}$ output for EDE (iff $\chi_{F2}^{(2)}>\chi_{F2}^{(1)}$)method:\\
$$
[j_1^{(1)},j_2^{(1)}],\,\chi_{F1}^{(1)}=x_{j_1^{(1)}},\,\chi_{F2}^{(1)}=x_{j_2^{(1)}},\,\chi_{D}^{(1)}=\frac{\chi_{F1}^{(1)}+\chi_{F2}^{(1)}}{2} \\
$$
\end{enumerate}
We can also try to produce output for ESE method for each BEDE iteration, but it is not always working, since $[\chi_{F1}^{(i)},\chi_{F2}^{(i)}]$ does not necessary contain $[x_r,x_l]$. 
Similarly we can produce output for EDE methods for each BESE iteration but again we are not sure that $[x_r^{(i)},x_l^{(i)}]\supset[\chi_{F1}^{(i)},\chi_{F2}^{(i)}]$.

\section{Experiments and Results}
We design small experiments by taking a suitable smooth function of known inflection point p, an interval $[a,b]$ that covers all the possible cases $p\in[a,b],\,p<a,\,p>b$ and we add a uniform error $\epsilon_i\sim{U(-r,r)}$ via the process \ref{eq:errdat}. 
\subsection{Symmetric sigmoid curves}
We find the points $x_1,x_{99}$ that give the first $1\%$ and the $99\%$ of sigmoid' s capacity L because those points have economic sense. So our interval $[a,b]$ is always relevant to the interval $[x_1,x_{99}]$. 
\subsubsection{The Fisher-Pry sigmoid curve with total symmetry}
\begin{figure}[hbtp]
\begin{center}
\includegraphics[scale=0.5,angle=0]{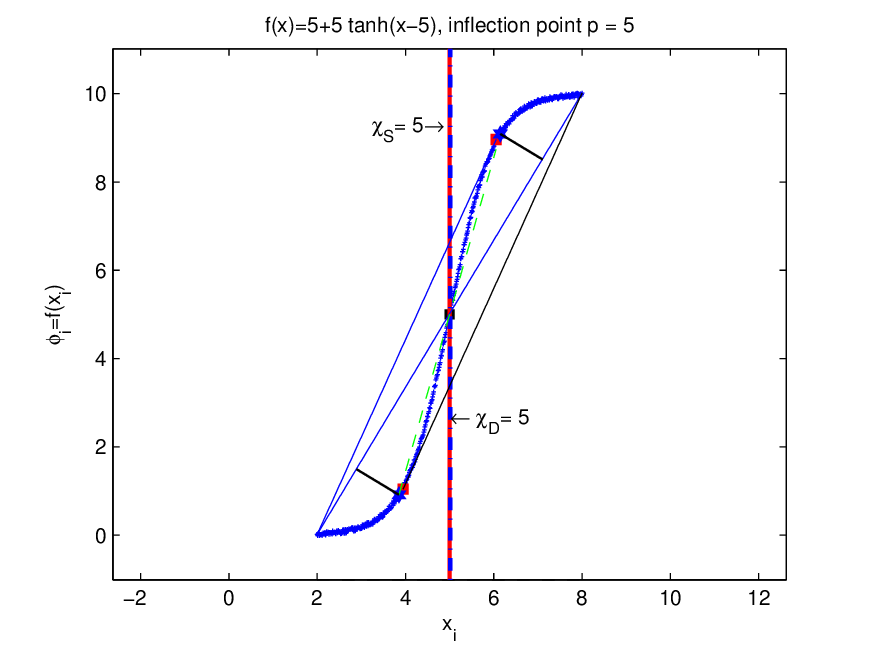}
\caption{(color online, scaled) Fisher-Pry sigmoid with total symmetry and without error}\label{fig:fish51no}
\end{center}
\end{figure}

Let' s take the function:
\beq
f \left( x \right) =5+5\,\tanh \left( x-5 \right)
\eeq
after \cite{fis-71}, which has $p=5,\,L=10,\,x_{1}=2.7024,\,x_{99}=7.2976$ and examine it at the interval $[2,8]$ in order to have data symmetry w.r.t. inflection point. The function is also symmetrical around inflection point, i.e. we have \emph{total symmetry}.\\
From Corollary \ref{cor:xtanlr} we compute $x_l=5.970315941,\,x_r=4.029684059$,\\
${x_{F}}_{1}=3.850750196,\,{x_{F}}_{2}=6.149249804$, all inside $[2,8]$, thus all methods are theoretically applicable.
We first take $n=500$ sub-intervals equal spaced without error just for checking our estimators. The results are presented at Table \ref{tab:fis-symnoer}.
\begin{table}[htbp]
\caption{Fisher-Pry sigmoid, p=5, total symmetry, n=500, no-error}
\label{tab:fis-symnoer}
{\begin{tabular}{ccccc} \toprule
$j_{{r}}=170$&$j_{{l}}=332$&$\chi_{{r}}= 4.028$&$\chi_{{l}}= 5.972$&$\chi_{{S}}= 5.0$\\ 
$j_{F1}=155$&$j_{F2}=347$&$\chi_{{F1}}= 3.884$&$\chi_{{F2}}= 6.1520$&$\chi_{{D}}= 5.0$\\ 
\botrule
\end{tabular}}
\end{table}
We observe that $\chi_l=5.9720,\,\chi_r=4.0280,\,{\chi_{F}}_{1}=3.8480,\,{\chi_{F}}_{2}=6.1520$ are very close to the theoretically expected values, so we are on the results of Lemma \ref{lem:esecon}. The absolutely accuracy from the first apply of all methods confirms our theoretical analysis. All important lines, curves and points are presented at Figure \ref{fig:fish51no}. 

\begin{figure}[hbtp]
\begin{center}
\includegraphics[scale=0.5,angle=0]{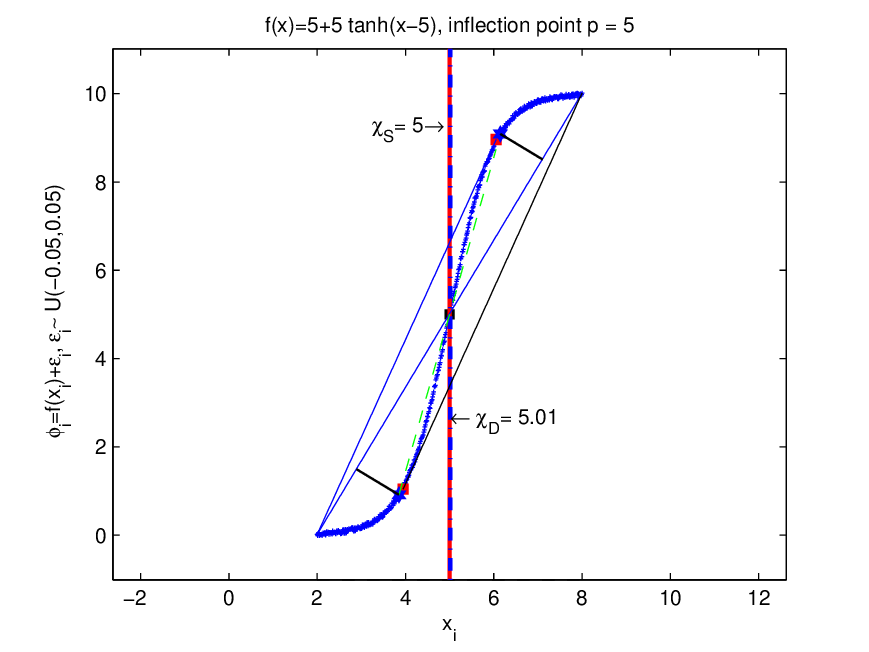}
\caption{(color online, scaled) Fisher-Pry sigmoid with total symmetry and error r=0.05}\label{fig:fish51}
\end{center}
\end{figure}

We next add the error term $\epsilon_i\sim{U(-0.05,0.05)}$ via the process \ref{eq:errdat} and run our algorithms again.The results are presented at Table \ref{tab:fis-symer}.   \\ 
\begin{table}[htbp]
\caption{Fisher-Pry sigmoid, p=5, data symmetry, $[a,b]=[2,8]$, n=500, error r=0.05}
\label{tab:fis-symer}
{\begin{tabular}{ccccc} \toprule
$j_{{r}}=163$&$j_{{l}}=339$&$\chi_{{r}}=3.944$&$\chi_{{l}}=6.056$&$\chi_{{S}}=5.000$\\ 
$j_{F1}=155$&$j_{F2}=349$&$\chi_{{F1}}=3.848$&$\chi_{{F2}}=6.176$&$\chi_{{D}}=5.012$\\ 
\botrule
\end{tabular}}
\end{table}
Again the estimations are close to the theoretically expected and both methods gave the true answer from the first apply. We present the ESE and EDE intervals and estimators together with the true function and the noisy data at Figure \ref{fig:fish51},
where we present all important points $x_l,x_r,x_{{F1}},x_{{F2}}$ with the relevant tangent lines and the size of the minimum/maximum of $F$.
Due to the \emph{total symmetry} the (dashed) line connecting $(x_{{F1}},f(x_{{F1}})),(x_{{F2}},f(x_{{F2}}))$ passes
from the point $(p,f(p))$.

\subsubsection{The Fisher-Pry sigmoid curve with data left asymmetry}
We continue with the same sigmoid function, but now we form proper our $[a,b]$ to show data asymmetry w.r.t. inflection point.\\
Let' s take for example $[4.2,8]$. If we do our theoretical computations we find $x_l=5.974322740,\,x_r=4.029684059$, ${x_{F}}_{1}=4.025677260,\,{x_{F}}_{2}=5.974322740$. We have that $x_r<a$, so $\chi_r$ has to estimate $a=4.2$ and $\chi_S$ must be close to $4.703504993$. Additionally, ${x_{F}}_{1}<a$, so ${\chi_{F}}_{1}$ must be also an estimation of $a$, thus $\chi_D$ must lie near the value $5.087161370$. It' s time to see if our theoretical predictions will be confirmed by experiment.\\
We use for comparability the same Standard Partition as before and have the output presented at Table \ref{tab:fis-leftasymno}. We have confirmation of our theory. 
\begin{table}[htbp]
\caption{Fisher-Pry sigmoid, p=5, data left symmetry, $[a,b]=[4.2,8]$, n=500, no error}
\label{tab:fis-leftasymno}
{\begin{tabular}{ccccc} \toprule
$j_{{r}}=2$&$j_{{l}}=156$&$\chi_{{r}}= 4.2076$&$\chi_{{l}}= 5.3780$&$\chi_{{S}}= 4.7928$\\ 
$j_{F1}=1$&$j_{F2}=234$&$\chi_{{F1}}= 4.2$&$\chi_{{F2}}= 5.9708$&$\chi_{{D}}= 5.0854$\\ 
\botrule
\end{tabular}}
\end{table}

\par It is time now to try iterations based on ESE and EDE intervals that contain inflection point and to observe remarkable convergence to the real value of $p=5$ for both methods. We present ESE \& EDE iterations at Table \ref{tab:fisno-iters}.

\begin{table}[htbp]
\caption{ESE \& EDE iterations for Fisher-Pry sigmoid, p=5, data left asymmetry, $[a,b]=[4.2,8]$, n=500, no-error}
\label{tab:fisno-iters}
{\begin{tabular}{ccccccc} \toprule
&(a) ESE&&&&(b) EDE&\\ %\colrule
$\chi_{{r}}$&$\chi_{{l}}$&$\chi_{{S}}$& &$\chi_{F1}$&$\chi_{F2}$&$\chi_{{D}}$\\ \colrule
4.8156& 5.3704& 5.0930& &4.5192& 5.4844& 5.0018\\ 
4.8232& 5.0892& 4.9562& &4.7244& 5.2716& 4.9980\\ 
4.9524& 5.0892& 5.0208& &4.8460& 5.1576& 5.0018\\ 
4.9600& 5.0208& 4.9904& &4.9068& 5.0892& 4.99800\\  
4.9904& 5.0208& 5.0056& &4.9448& 5.0512& 4.99800\\ \botrule
\end{tabular}}
\end{table}

\par Let' s add the same error term $\epsilon_i\sim{U(-0.05,0.05)}$ and run our algorithms. The results at Table \ref{tab:fis-leftasymer} clearly are close enough to the theoretical expectations.\\
\begin{table}[htbp]
\caption{Fisher-Pry sigmoid, p=5, data left symmetry, $[a,b]=[4.2,8]$, n=500, error r=0.05}
\label{tab:fis-leftasymer}
{\begin{tabular}{ccccc} \toprule
$j_{{r}}=3$&$j_{{l}}=149$&$\chi_{{r}}=4.2304$&$\chi_{{l}}=5.3248$&$\chi_{S}=4.7700$\\
$j_{F1}=3$&$j_{F2}=231$&$\chi_{{F1}}=4.2152$&$\chi_{{F2}}=5.9480$&$\chi_{{D}}=5.0816$\\ \botrule
\end{tabular}}
\end{table} 

We observe that ESE method did not estimate the inflection point with acceptable accuracy, so it is time to run the ESE and EDE iterations. The results, Table \ref{tab:fis-leftasymer-iters}  show a clear improvement of both estimations. 

\begin{table}[htbp]
\caption{ESE \& EDE iterations for Fisher-Pry sigmoid, p=5, data left asymmetry, $[a,b]=[4.2,8]$, n=500, error r=0.05}
\label{tab:fis-leftasymer-iters}
{\begin{tabular}{ccccccc} \toprule
&(a) ESE&&&&(b) EDE&\\ %\colrule
$\chi_{{r}}$&$\chi_{{l}}$&$\chi_{{S}}$& &$\chi_{F1}$&$\chi_{F2}$&$\chi_{{D}}$\\ \colrule
4.8156&5.3248& 5.0702& &4.5268& 5.5148& 5.0208\\ 
4.9144&5.1576& 5.0360& &4.7244& 5.2412& 4.9828\\ 
\botrule
\end{tabular}}
\end{table} 

\begin{figure}[hbtp]
\begin{center}
\includegraphics[scale=0.5,angle=0]{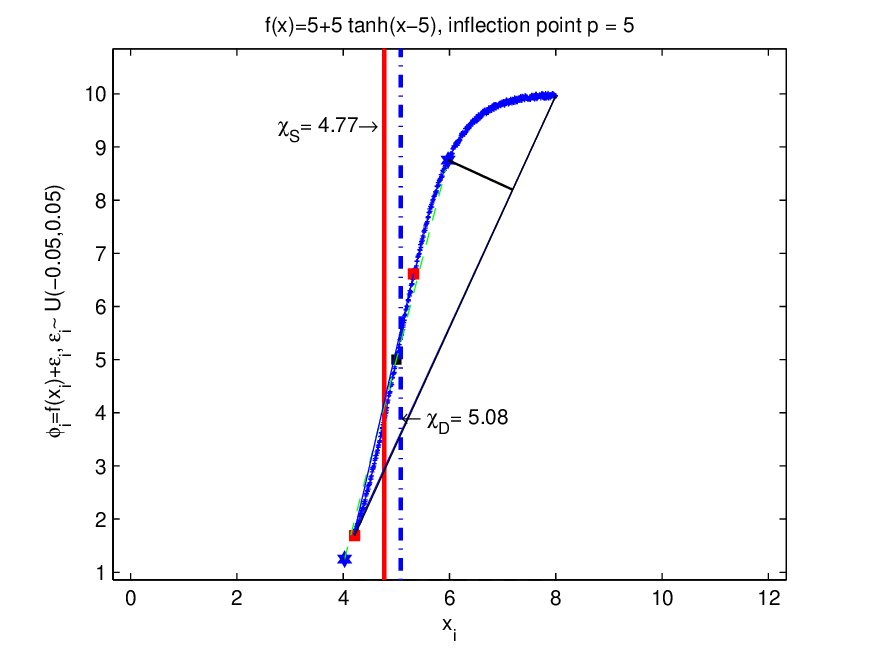}
\caption{(color online, scaled) Fisher-Pry sigmoid with data left asymmetry and error r=0.05}\label{fig:fish52}
\end{center}
\end{figure}

All the points of interesting are presented at Figure \ref{fig:fish52}, where we see that interval does not contain both $x_l,x_r$ and $x_{{F1}},x_{{F2}}$.

%\clearpage

\subsection{Non symmetric sigmoid curves}
\par We continue our study with non symmetric sigmoid curves that appear in Economics and other disciplines.

\subsubsection{The Gompertz sigmoid curve}
Let' s examine the function:
\beq
f \left( x \right) =10\,{{\rm e}^{-{{\rm e}^{5}}{{\rm e}^{-x}}}}
\eeq
after \cite{gom-25}, in the interval $[3.5,8]$. The basic properties are presented at Table \ref{tab:gom-prop}.\\
\begin{table}[htbp]
\caption{Basic properties of a Gompertz sigmoid}
\label{tab:gom-prop}
{\begin{tabular}{ccc} \toprule
$L=10$&$x_1=3.472820374$&$x_{99}=9.600149227$\\ 
$x_{{r}}=4.138928270$&$x_{{l}}=5.887451706$&$x_{{S}}=5.013189988$\\ 
$x_{{F1}}=4.095750735$&$x_{{F2}}=6.290768183$&$x_{{D}}=5.193259460$\\ 
\botrule
\end{tabular}}
\end{table}
It is easy to prove that f is $(0.224,1.0)$-asymptotically symmetric around inflection point, so we can handle it similar to a symmetric sigmoid only for a distance of $\pm{1}$ from $p=5$.
\par We use, for comparison reasons, the same SP with 500 sub-intervals without error and obtain the Table \ref{tab:gom-noer} which is absolutely compatible with theoretical predictions. The ESE \& EDE iterations are showed at Table \ref{tab:gom-noer-iters} where we observe convergence to the real p for both two methods. 
\begin{table}[htbp]
\caption{Gompertz sigmoid, p=5, asymmetry, n=500, no-error}
\label{tab:gom-noer}
{\begin{tabular}{ccccc} \toprule
$j_{{r}}=72$&$j_{{l}}=266$&$\chi_{{r}}=4.139$&$\chi_{{l}}=5.885$&$\chi_{{S}}=5.012$\\ 
$j_{F1}=67$&$j_{F2}=311$&$\chi_{{F1}}=4.094$&$\chi_{{F2}}=6.290$&$\chi_{{D}}=5.192$\\ \botrule
\end{tabular}}
\end{table}

\begin{figure}[hbtp]
\begin{center}
\includegraphics[scale=0.5,angle=0]{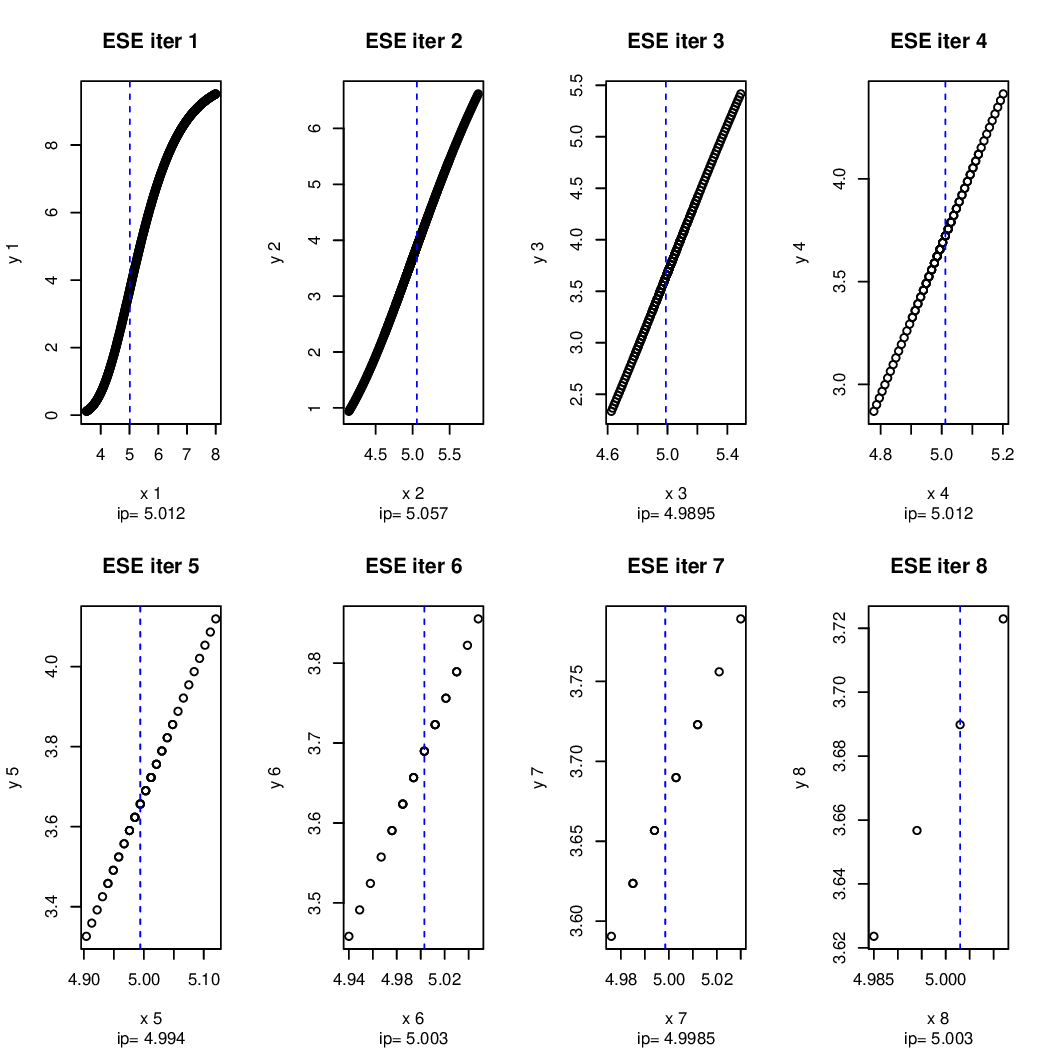}
\caption{ESE iterations Gompertz asymmetrical sigmoid no error}\label{fig:eseiter}
\end{center}
\end{figure}

\begin{figure}[hbtp]
\begin{center}
\includegraphics[scale=0.5,angle=0]{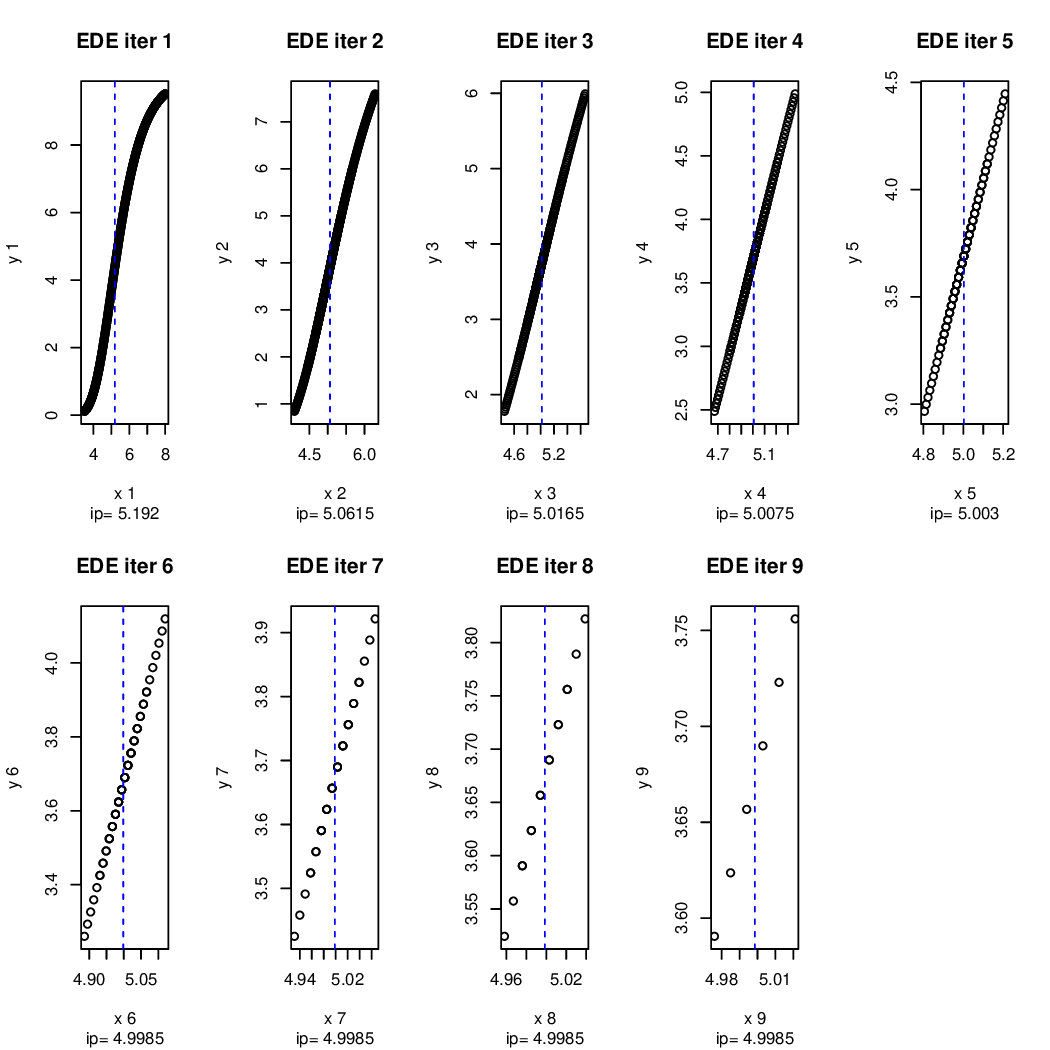}
\caption{EDE iterations Gompertz asymmetrical sigmoid no error}\label{fig:edeiter}
\end{center}
\end{figure}

\begin{table}[htbp]
\caption{ESE \& EDE iterations for Gompertz sigmoid, p=5, asymmetry, $[a,b]=[3.5,8]$, n=500, no error}
\label{tab:gom-noer-iters}
{\begin{tabular}{ccccccc} \toprule
&(a) ESE&&&&(b) EDE&\\ %\colrule
$\chi_{{r}}$&$\chi_{{l}}$&$\chi_{{S}}$& &$\chi_{F1}$&$\chi_{F2}$&$\chi_{{D}}$\\ \colrule
4.625& 5.489&5.0570& &4.4540& 5.6690& 5.0615\\ 
4.778& 5.201&4.9895& &4.6700& 5.3630& 5.0165\\  
4.904& 5.120&5.0120& &4.8050& 5.2100& 5.0075\\ 
4.940 &5.048&4.9940& &4.8860& 5.1200& 5.0030\\
4.976 &5.030&5.0030& &4.9310& 5.0660& 4.9985\\ 
4.985 &5.012&4.9985& &4.9580& 5.0390& 4.9985\\
\botrule
\end{tabular}}
\end{table} 
We can watch the iteration convergence of the two methods in the non error case to the actual value of inflection point at Figure \ref{fig:eseiter} and Figure \ref{fig:edeiter}, where we have left the methods to stop when the minimum required points are reached.

%\clearpage

\begin{figure}[hbtp]
\begin{center}
\includegraphics[scale=0.5,angle=0]{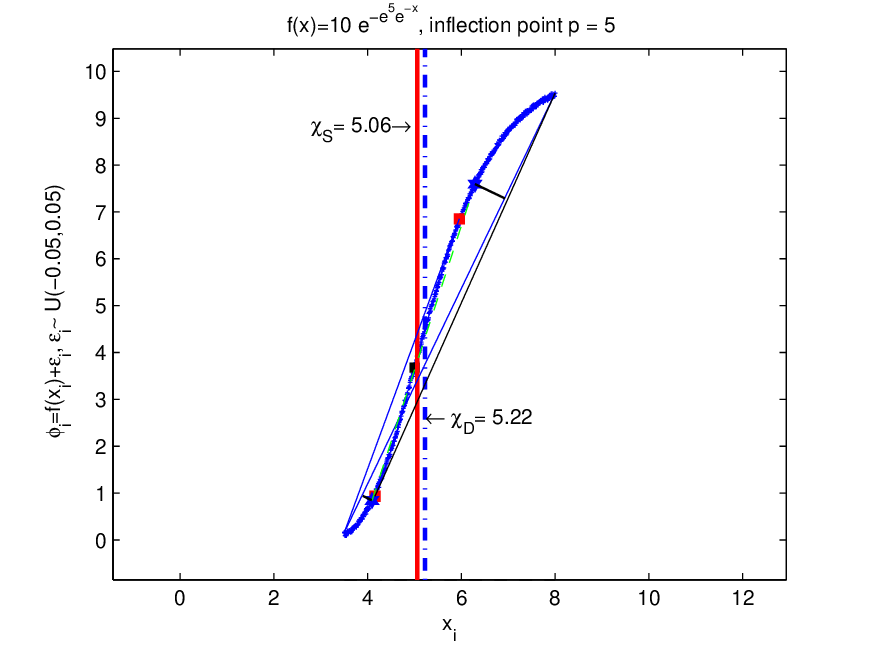}
\caption{(color online, scaled) Gompertz asymmetrical sigmoid with error r=0.05}\label{fig:gom51}
\end{center}
\end{figure}

We continue with our familiar SP by adding error uniformly distributed by $U(-0.05,0.05)$ and the results are given at Table \ref{tab:gom-er} while ESE \& EDE iterations are shown at Table \ref{tab:gom-er-iters}. From these Tables we conclude that convergence to the true value of inflection point $p=5$ occurs from the iterative application of ESE and EDE methods in one or two steps only. All points of interest and data are presented at Figure \ref{fig:gom51}.\\  
\begin{table}[htbp]
\caption{Gompertz sigmoid, p=5, asymmetry, $[a,b]=[3.5,8]$, n=500, error r=0.05}
\label{tab:gom-er}
{\begin{tabular}{ccccc} \toprule
$j_{{r}}=74$&$j_{{l}}=274$&$\chi_{{r}}=4.1570$&$\chi_{{l}}=5.9570$&$\chi_{{S}}=5.0570$\\ 
$j_{F1}=66$&$j_{F2}=319$&$\chi_{{F1}}=4.0850$&$\chi_{{F2}}=6.3620$&$\chi_{{D}}=5.2235$\\ \botrule
\end{tabular}}
\end{table}

\begin{table}[htbp]
\caption{ESE \& EDE iterations for Gompertz sigmoid, p=5, asymmetry, $[a,b]=[3.5,8]$, n=500, error r=0.05}
\label{tab:gom-er-iters}
{\begin{tabular}{ccccccc} \toprule
&(a) ESE&&&&(b) EDE&\\ %\colrule
$\chi_{{r}}$&$\chi_{{l}}$&$\chi_{{S}}$& &$\chi_{F1}$&$\chi_{F2}$&$\chi_{{D}}$\\ \colrule
4.6340& 5.5340& 5.0840& &4.472& 5.642& 5.057\\ 
4.8590& 5.1560& 5.0075& & & & \\  
\botrule
\end{tabular}}
\end{table} 

%\clearpage

\subsection{Non sigmoid curves}
Our analysis is applicable also to non sigmoid curves, not necessary symmetric or with data symmetry. We shall proceed by making two experiments for a symmetric $3^{rd}$ order polynomial.

\subsubsection{A symmetric $3^{rd}$ order polynomial with total symmetry}
Let the polynomial function:
\beq \label{eq:pol3sym}
f \left( x \right) =-\frac{1}{3}\,{x}^{3}+\frac{5}{2}\,{x}^{2}-4\,x+\frac{1}{2}
\eeq
We study it at $[-2,7]$, it has inflection point at $p=2.5$ and we have \emph{total symmetry} while the interesting points are presented at Table \ref{tab:3rd-prop}.\\
\begin{table}[htbp]
\caption{Basic properties of a $3^{rd}$ order polynomial with total symmetry}
\label{tab:3rd-prop}
{\begin{tabular}{lll} \toprule
$x_{{r}}=0.25$&$x_{{l}}=4.75$&$x_{{S}}=2.50$\\ 
$x_{{F1}}=-0.09807621078$&$x_{{F2}}=5.098076211$&$x_{{D}}=2.50$\\ \botrule
\end{tabular}}
\end{table}

\par The SP with 500 sub-intervals without error gives Table \ref{tab:pol3-noer} which is absolutely compatible with theoretical predictions. There is no need for any kind of iteration, because both methods agree with the true value.

\begin{table}[htbp]
\caption{$3^{rd}$ order polynomial, total symmetry, p=2.5, n=500, no-error}
{\begin{tabular}{ccccc} \toprule
$j_{{r}}=126$&$j_{{l}}=376$&$\chi_{{r}}=0.25$&$\chi_{{l}}=4.75$&$\chi_{{S}}=2.50$\\ 
$j_{F1}=107$&$j_{F2}=395$&$\chi_{{F1}}=-0.092$&$\chi_{{F2}}=5.092$&$\chi_{{D}}=2.50$\\ \botrule
\end{tabular}}
\label{tab:pol3-noer}
\end{table}

\par
\begin{figure}[hbtp]
\begin{center}
\includegraphics[scale=0.5,angle=0]{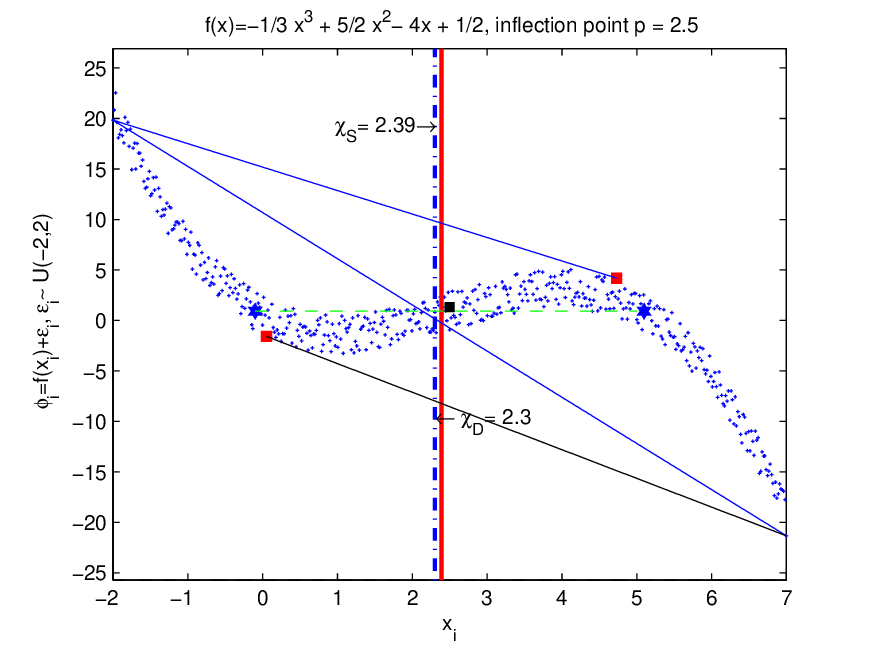}
\caption{(color online, unscaled) Symmetric $3^{rd}$ order polynomial with error r=2.0}\label{fig:pol351}
\end{center}
\end{figure}

 The same SP with uniform error distributed by $U(-2,2)$ gives the results of Table \ref{tab:pol3-er} and two ESE iterations are presented at Table \ref{tab:pol3eseit}. All points and data are shown at Figure \ref{fig:pol351}. 

\begin{table}[htbp]
\caption{Symmetric $3^{rd}$ order polynomial, total symmetry, p=2.5, n=500, error r=2.0}
\label{tab:pol3-er}
{\begin{tabular}{ccccc} \toprule
$j_{{r}}=115$&$j_{{l}}=375$&$\chi_{{r}}=0.052$&$\chi_{{l}}=4.732$&$\chi_{{S}}=2.392$\\ 
$j_{F1}=105$&$j_{F2}=375$&$\chi_{{F1}}=-0.128$&$\chi_{{F2}}=4.732$&$\chi_{{D}}=2.302$\\ \botrule
\end{tabular}}
\end{table}

\begin{table}[htbp]
\caption{ESE iterations for $3^{rd}$ order polynomial, p=2.5, total symmetry, n=500, error r=2.0}
\label{tab:pol3eseit}
{\begin{tabular}{ccc} \toprule
$\chi_{{r}}$&$\chi_{{l}}$&$\chi_{{S}}$\\ 
1.222& 3.688& 2.455\\ 
1.564& 3.382& 2.473\\ \botrule
\end{tabular}}
\end{table}

\subsubsection{A symmetric $3^{rd}$ order polynomial with data right asymmetry}
For the same symmetric $3^{rd}$ order polynomial \ref{eq:pol3sym} we change the interval to $[-2,8]$, thus we have data right asymmetry now. Our critical points are written here at Table \ref{tab:3rd-drasym}.

\begin{table}[htbp]
\caption{Basic properties of a $3^{rd}$ order polynomial with data right asymmetry}
\label{tab:3rd-drasym}
{\begin{tabular}{lll} \toprule
$x_{{r}}=-0.25$&$x_{{l}}=4.75$&$x_{{S}}=2.25$\\ 
$x_{{F1}}=-0.429732639$&$x_{{F2}}=5.429732639$&$x_{{D}}=2.50$\\ \botrule
\end{tabular}}
\end{table}

The case of SP with 500 sub-intervals and no error gives Table \ref{tab:pol3-noer2}, while ESE and EDE iterations are presented at Table \ref{tab:pol3-no2-iters}.  First results are absolutely compatible with theoretical predictions for ESE method. %For example we are waiting that $\chi_S=2.25$ and we found $2.24$. 

\begin{table}[htbp]\small
\caption{Symmetric $3^{rd}$ order polynomial, data right asymmetry, p=2.5, n=500, $[-2,8]$, no-error}
\label{tab:pol3-noer2}
{\begin{tabular}{ccccc} \toprule
$j_{{r}}=89$&$j_{{l}}=338$&$\chi_{{r}}=- 0.24$&$\chi_{{l}}=4.74$&$\chi_{{S}}=2.25$\\ 
$j_{F1}=80$&$j_{F2}=372$&$\chi_{{F1}}=- 0.42$&$\chi_{{F2}}=5.42$&$\chi_{{D}}=2.50$\\ \botrule
\end{tabular}}
\end{table}

\begin{table}[htbp]
\caption{ESE \& EDE iterations for $3^{rd}$ order polynomial, p=5, p=2.5, n=500, $[-2,8]$, no-error}
\label{tab:pol3-no2-iters}
{\begin{tabular}{ccccccc} \toprule
&(a) ESE&&&&(b) EDE&\\ %\colrule
$\chi_{{r}}$&$\chi_{{l}}$&$\chi_{{S}}$& &$\chi_{F1}$&$\chi_{F2}$&$\chi_{{D}}$\\ \colrule
1.3800& 3.8800& 2.6300& &0.82000& 4.1800& 2.5000\\ 
1.8000& 3.0600& 2.4300& & & & \\  
2.2200& 2.8600& 2.5400& & & & \\
2.3200& 2.6400& 2.4800& & & & \\
2.4200& 2.5800& 2.5000& & & & \\
2.4600& 2.5400& 2.5000& & & & \\
\botrule
\end{tabular}}
\end{table} 

\par
\begin{figure}[hbtp]
\begin{center}
\includegraphics[scale=0.5,angle=0]{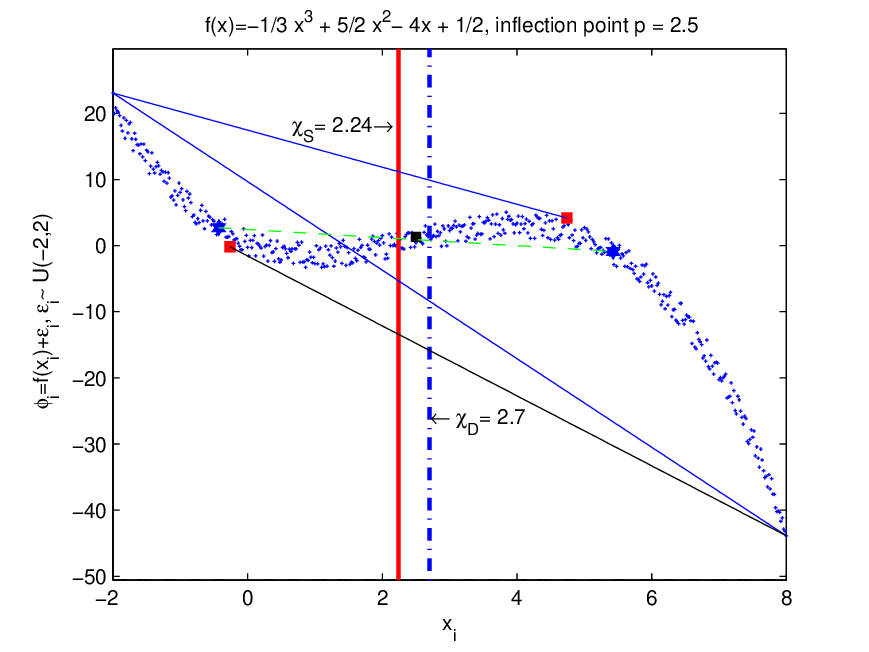}
\caption{(color online, unscaled) Symmetric $3^{rd}$ order polynomial with data right asymmetry}\label{fig:pol352}
\end{center}
\end{figure}

 We add uniform error distributed by $U(-2,2)$ and we have the results of Table \ref{tab:pol3-er2}, while one ESE \& one EDE iteration are given at Table \ref{tab:pol3-er2-iters}.
 All points and data are presented at Figure \ref{fig:pol352}. 

\begin{table}[htbp]\small
\caption{Symmetric $3^{rd}$ order polynomial, data right asymmetry, p=2.5, n=500, $[-2,8]$, error r=2.0}
\label{tab:pol3-er2}
{\begin{tabular}{ccccc} \toprule
$j_{{r}}=88$&$j_{{l}}=338$&$\chi_{{r}}=- 0.26$&$\chi_{{l}}=4.74$&$\chi_{{S}}=2.24$\\
$j_{F1}=88$&$j_{F2}=384$&$\chi_{{F1}}=- 0.26$&$\chi_{{F2}}=5.66$&$\chi_{{D}}=2.70$\\ \botrule
\end{tabular}}
\end{table}

\begin{table}[htbp]
\caption{ESE \& EDE iterations for $3^{rd}$ order polynomial, p=2.5, n=500, $[-2,8]$, error r=2.0}
\label{tab:pol3-er2-iters}
{\begin{tabular}{ccccccc} \toprule
&(a) ESE&&&&(b) EDE&\\ %\colrule
$\chi_{{r}}$&$\chi_{{l}}$&$\chi_{{S}}$& &$\chi_{F1}$&$\chi_{F2}$&$\chi_{{D}}$\\ \colrule
1.46& 3.84& 2.65& &0.86& 3.84& 2.35\\ \botrule
\end{tabular}}
\end{table} 

\par There exist a problem here. Although we have a symmetric polynomial, the TESE is not equal to the true inflection point. A remedy for this problem for the class of $3^{rd}$ order polynomials is given with the next Lemma.\\ 

\begin{lemma}{The \rm $3^{rd}$ order polynomial ESE correction.}\\
Let a $3^{rd}$ order polynomial $f \left( x \right) =\alpha\,{x}^{3}+\beta\,{x}^{2}+\gamma\,x+\delta$ and let p its inflection point. Then it holds exactly that:
\beq
p=\frac{1}{3}\,x_{{l}}+\frac{1}{3}\,x_{{r}}+\frac{1}{6}\,a+\frac{1}{3}\,b
\eeq
and 
$$
\hat{p}=\frac{1}{3}\,\chi_{{l}}+\frac{1}{3}\,\chi_{{r}}+\frac{1}{6}\,a+\frac{1}{3}\,b
$$
is a consistent estimator of trapezoidal estimated p.
\end{lemma}
\textbf{Proof}\\
The inflection point because $\alpha\neq{0}$ is found from the root of the second derivative, i.e. $6\,\alpha\,p+2\,\beta=0$ or $p=-\,\frac {\beta}{3\,\alpha}$. Due to Corollary \ref{cor:xtanlr} we have for the $x_l$ that:
$$
3\,\alpha\,{x}^{2}+2\,\beta\,x+\gamma={\frac {\alpha\,{x}^{3}+\beta\,{
x}^{2}+\gamma\,x-\alpha\,{a}^{3}-\beta\,{a}^{2}-\gamma\,a}{x-a}}
$$  
or
$$
\left( x-a \right) ^{2} \left( \alpha\,a+\beta+2\,\alpha\,x \right) =0
$$
so the internal solution $x_l$ is:
\beq
x_l=-\frac {\alpha\,a+\beta}{2\,\alpha}
\eeq
For the $x_r$ we have similar that:
$$
3\,\alpha\,{x}^{2}+2\,\beta\,x+\gamma={\frac {\alpha\,{b}^{3}+\beta\,{
b}^{2}+\gamma\,b-\alpha\,{x}^{3}-\beta\,{x}^{2}-\gamma\,x}{b-x}}
$$  
or
$$
\left( b-x \right) ^{2} \left( \alpha\,b+\beta+2\,\alpha\,x \right)=0
$$
so the internal solution $x_r$ is:
\beq
x_r=-\frac {\alpha\,b+\beta}{2\,\alpha}
\eeq
By adding $x_l$ and $x_r$ we obtain:
$$
x_{{l}}+x_{{r}}=-\frac{1}{2}\,a-{\frac {\beta}{\alpha}}-\frac{1}{2}\,b
$$
and if we remember that $p=-\,\frac {\beta}{3\,\alpha}$ we obtain:
$$
x_{{l}}+x_{{r}}=-\frac{1}{2}\,a+3\,p-\frac{1}{2}\,b
$$
or finally
\beq
p=\frac{1}{3}\,x_{{l}}+\frac{1}{3}\,x_{{r}}+\frac{1}{6}\,a+\frac{1}{3}\,b
\eeq
\par Since we have proven that $\chi_l,\chi_r$ are consistent estimators of trapezoidal calculated values of $x_l,x_r$ we can take a consistent estimation for trapezoidal calculated p by replacing the unknown $x_l,x_r$ with the estimators $\chi_l,\chi_r$. \hfill\(\Box\) \\

\par As an example, we come back to the case of $3^{rd}$ order symmetric polynomial with data right asymmetry. We have that $a=-2,\,b=8$ and from Table \ref{tab:pol3-er2} is $\chi_{{r}}=- 0.26,\,\chi_{{l}}=4.74$, so we have that:
$$
\hat{p}={\frac{1}{3}\,\chi_{{l}}+\frac{1}{3}\,\chi_{{r}}+\frac{1}{6}\,a+\frac{1}{3}\,b}=2.493333333
$$ 
which is much closer to the true value of $2.5$. 
\par The above analysis can be extended to every function, if we can find analytically a relation between inflection point and $x_l,x_r,a,b$.

\section{Discussion}
The sigmoid or S-shaped pattern is common in many disciplines: utility theory \cite{fri-48} \& technological substitution models \cite{mar-79,mod-92} in Economics, growth theory \cite{ber-60} \& allosterism \cite{bar-75} in Biology, population dynamics \cite{ver-38,tur-02} in Ecology, titration data analysis \cite{moh-77} \& autocatalysis \cite{upa-06} in Analytical Chemistry, dose-response \cite{med-89} in Medicine and many others. 
 Starting from the problem of identifying the inflection point for the data $\left\{\left(x_i,\phi_i\right),\,i=0,1,\ldots,n\right\}$ we have created two geometrical methods in order to solve this problem: the Extremum Surface Estimator (ESE method) and the Extremum Distance Estimator (EDE method). We did not perform neither regression nor splines representation analysis. In addition, no concepts like discrete or digital curvature were defined or used. Instead we focused on finding two points where the true inflection point lies between and then we took as estimator the middle point, just like in bisection method. \\
\par The methods can be applied iteratively similarly to bisection method and converge to the true inflection point either from the first or after a few iterations only. The R Package {\it inflection} is available at \cite{inf-13} for using ESE \& EDE methods with R. Other implementations for the methods have been done with FORTRAN, Maple and Matlab. If we have large data sets then FORTRAN is more efficient, for example the problem of estimating the inflection point of $n=10000$ data pairs needed less than 7 sec CPU time in FORTRAN GNU compiler and a typical Intel Core i5 CPU with 4 GB RAM memory.

\end{document}